\documentclass[a4paper,12pt]{amsart}
\usepackage[utf8]{inputenc}
\usepackage[centering]{geometry}
\usepackage{amsmath, amsthm, amssymb, amsfonts} 
\usepackage{mathtools}
\mathtoolsset{showonlyrefs}
\usepackage{graphicx}
\usepackage[T1]{fontenc}
\usepackage{booktabs}
\usepackage{xcolor}
\usepackage{mathrsfs}

\usepackage{romannum}

\usepackage{enumitem}

\usepackage{hyperref}
\hypersetup{
    colorlinks=true,
    linkcolor=blue,
    filecolor=magenta,      
    urlcolor=cyan,
    pdftitle={Overleaf Example},
    pdfpagemode=FullScreen,
    unicode=true
}
\theoremstyle{plain}

\theoremstyle{thmstyleone}%
\newtheorem{theorem}{Theorem}[section]
\newtheorem{proposition}[theorem]{Proposition}%
\newtheorem{example}[theorem]{Example}%
\newtheorem{remark}[theorem]{Remark}%
\newtheorem{lemma}[theorem]{Lemma}
\newtheorem{definition}[theorem]{Definition}
\newtheorem{corollary}[theorem]{Corollary}
\newtheorem{claim}[theorem]{Claim}
\newtheorem{conjecture}{Conjecture}%
\newtheorem{maintheorem}{Theorem}%

\def \FF {{\mathbb F}}

\title{Algebraic Exceptional Set of a Three-Component Curve on Hirzebruch Surfaces}
\author{Wei Chen}

\address{ Department of Mathematics and Physics \\
Roma Tre University   \\
Largo San Leonardo Murialdo, I-00146 \\
Rome, Italy}
\email{wei.chen@uniroma3.it, weichen97.ag@gmail.com}

\begin{document}

\pagenumbering{arabic}

\begin{abstract}
We study the algebraic exceptional set of a three-component curve $B$ with normal crossings on a Hirzebruch surface $\mathbb{F}_e$. If $K_{\mathbb{F}_{e}}+B$ is big and no component of $B$ is a fiber or the rational curve with negative self-intersection, we prove that the algebraic exceptional set is finite, and in most cases give it an effective bound. We also prove that the algebraic exceptional set
coincides with the set of curves that are hyper-bitangent to $B$. 
\end{abstract}

\maketitle

\tableofcontents 

\section{Introduction}
\label{sec:intro}
Rational curves on smooth algebraic surfaces are fundamental geometric objects in the study of algebraic surfaces. It is conjectured that there are only finitely many rational curves on a complex smooth projective surface $S$ of general type. In this case, Bogomolov's result \cite{Bogomolov77}  implies that there are only finitely many rational curves on $S$ if $c_{1}^{2}(S)-c_{2}(S)>0$, and later Lu and Miyaoka \cite{Lu1995BoundingCI} showed that there are only finitely many smooth rational curves on any surface $S$ of general type. 

Similarly, rational curves on log surfaces are also of great importance. We are particularly interested in the \textit{algebraic exceptional set}, first introduced by Lang \cite{Lang86}, which is defined as follows:
\begin{definition}
    Let $(S, B)$ be a pair, where $S$ is a complex smooth projective surface and $B$ is a reduced curve on $S$. The algebraic exceptional set associated with this pair is defined as 
    \begin{equation*}
        \mathcal{E}(S, B) = 
        \left\{
            C \subset S 
            \ \vert \
            C \text{ is an integral rational curve such that }
                \#\nu_{C}^{-1}(C \cap B) \leq 2
        \right\},
    \end{equation*}
    where $\nu_{C}\colon C^{\nu} \to C$ is the normalization map of $C$. When the surface $S$ is clear from the context, we denote this set simply as $\mathcal{E}(B)$.
\end{definition}

In particular, when the pair is of log-general type---equivalently, when $K_{S} + B$ is a big divisor on $S$, and $B$ is a simple normal crossing divisor, we have the following conjecture:
\begin{conjecture}\label{conjecture 1}
    If $(S, B)$ is a log-smooth surface of log-general type, then $\mathcal{E}(B)$ is a finite set.
\end{conjecture}

From an arithmetic point of view, Conj.~\ref{conjecture 1} can be deduced from Vojta's conjecture in Diophantine geometry \cite[Conj.3.4.3]{vojta1987diophantine}; from an analytic-geometric perspective, it is related to the Green--Griffiths--Lang conjecture. For more details, we refer interested readers to \cite{Ascher2020}. It is worth mentioning that McQuillan \cite{PMIHES_1998__87__121_0} proved the Green-Griffiths-Lang conjecture for a complex smooth projective surface $S$ of general type with $c_{1}^{2}(S) - c_{2}(S) > 0$.

For the case when $S$ is the projective plane $\mathbb{P}^{2}$ and $B$ has normal crossing singularities at each intersection point of its irreducible components, it turns out that the fewer irreducible components $B$ has, the more challenging the problem becomes. If $B$ has at most two components, Conj.~\ref{conjecture 1} is still open, and the only known results are obtained assuming that the curve $B$ is very general (see \cite{XiChen04}, \cite{PacienzaRousseau+2007+221+235}, \cite{chen2023algebraichyperbolicitycomplementsgeneric}, \cite{ascher2024algebraicgreengriffithslangconjecturecomplements}); if $B$ has at least four components, Conj.~\ref{conjecture 1} is easy (see \cite[Prop.~4.3.1]{caporaso2024hypertangencyplanecurvesalgebraic}); if $B$ has three components, Conj.~\ref{conjecture 1} is true but requires a non-trivial proof (see \cite{CZ08}, \cite{guo2023vojtasabcconjecturealgebraic}). Recently, Caporaso and Turchet \cite{caporaso2024hypertangencyplanecurvesalgebraic} provided a new proof of the three-component case dropping the rationality assumption. Instead of $\mathcal{E}(B)$, they considered the set of integral curves that are \textit{hyper-bitangent} to $B$:
$$
    \operatorname{Hyp}(\mathbb{P}^{2},B,2) = 
    \left\{
        C \subset \mathbb{P}^{2} 
        \ \big\vert \ 
            C \text{ is an integral plane curve such that } \#\nu_{C}^{-1}(C \cap B) \leq 2 
    \right\},
$$
and they showed that $\operatorname{Hyp}(\mathbb{P}^{2},B,2) = \mathcal{E}(\mathbb{P}^2 , B)$. Moreover, they also proved that $\mathcal{E}(\mathbb{P}^2 , B)$ is finite and provided an effective bound for its cardinality.

Similarly, we define $\operatorname{Hyp}(S,B,2)$ for any pair $(S,B)$ as follows:
\begin{definition}
    Let $(S, B)$ be a pair, where $S$ is a complex smooth projective surface and $B$ is a reduced curve on $S$. The set of hyper-bitangent curves associated with this pair is defined as 
    \begin{equation*}
        \operatorname{Hyp}(S,B,2) = 
        \left\{
            C \subset S 
            \ \vert \ 
            C \text{ is an integral curve such that }
                 \#\nu_{C}^{-1}(C \cap B) \leq 2
        \right\}.
    \end{equation*}
    When the surface $S$ is clear from the context, we denote this set simply as $\operatorname{Hyp}(B,2)$.
\end{definition}

In this paper, we study the case when $S$ is a Hirzebruch surface $\mathbb{F}_{e},e\geq 0$ and $B$ is a curve with three irreducible components---in this case Conj.~\ref{conjecture 1} was open. Following Caporaso--Turchet's work\cite{caporaso2024hypertangencyplanecurvesalgebraic}, we adopt a slightly weaker assumption of the singularities on $B$ than simple normal crossing:
\begin{definition}\label{def : 3C-curve}
    Let $S$ be a smooth projective surface and $B$ be a curve on it. We say $B$ is a \textbf{$3C$-curve} if $B$ consists of three irreducible components such that it has normal crossing singularities at each intersection point of its components.
\end{definition}

We generalize and apply Caporaso--Turchet method \cite{caporaso2024hypertangencyplanecurvesalgebraic} and our results can be summarized in the following theorem:
\begin{maintheorem}\label{thm : main}
    Let $S$ be a Hirzebruch surface and $B$ a $3C$-curve on $S$ such that $K_{S}+B$ is a big divisor. If none of the irreducible components of $B$ is a fiber or the rational curve with negative self-intersection, then
    \begin{enumerate}
        \item $\mathcal{E}(B)=\operatorname{Hyp}(B,2)$;
        \item $\mathcal{E}(B)$ is a finite set.
    \end{enumerate}
\end{maintheorem}

This theorem is, to our knowledge, the first extension to a surface other than $\mathbb{P}^2$, of \cite{CZ08}'s finiteness result for the algebraic exceptional set of a three-component curve in $\mathbb{P}^2$. We also prove, under some extra assumptions, that 
$\mathcal{E}(B) = \emptyset$
if $B$ is general. See Prop.~\ref{prop: emptyness of exceptional set F0}, Prop.~\ref{prop: emptyness of exceptional set Fe}, and Prop.~\ref{prop: emptyness of exceptional set F1}.

The first part of Thm.~\ref{thm : main} is proved in Prop.~\ref{prop: F0 Hyp(B,2)=E(B)}, Prop.~\ref{prop: Fe Hyp(B,2)=E(B)}, and Prop.~\ref{prop: F1 Hyp(B,2)=E(B)}; the second part is established in Thm.~\ref{thm: F0 bound E(B)}, Thm.~\ref{thm: Fe bound E(B)}, and Thm.~\ref{thm: F1 bound E(B)}. In fact, we can bound $\mathcal{E}(B)$ effectively except for one special case on $\mathbb{F}_{1}$, where we obtain the finiteness of $\mathcal{E}(B)$ by applying the method developed by Corvaja and Zannier in \cite{CZ13},  and this will be discussed in part \ref{thm: F1 bound E(B)_b} of Thm.~\ref{thm: F1 bound E(B)}. Additionally, in Rem.~\ref{rmk: without Corvaja--Zannier}, we note that on $ \mathbb{F}_{e},e \geq 2 $, the Corvaja--Zannier method in \cite{CZ13} does not apply, in fact our method is independent of \cite{CZ13}.

\subsection*{Organization of the Paper} 
In Section~\ref{sec:prelim}, we study the properties of unibranch points in Subsection~\ref{subsec:unibranch points}, where we introduce Thm.~\ref{thm: bounddeltainvariant} and Thm.~\ref{thm: Strong Triangle Inequality}, which serve as key tools for our subsequent analysis. We then review basic facts about Hirzebruch surfaces in Subsection~\ref{subsec:Hirzebruch surfaces}, particularly focusing on a criterion to determine whether an integral divisor is big on Hirzebruch surfaces, as stated in Prop.~\ref{prop: Criterion for bigness on Hirzebruch Surfaces}.

In Section~\ref{sec:Hyper-bitangent Curves on Minimal Hirzebruch Surfaces}, we investigate $\operatorname{Hyp}(\mathbb{F}_{e},B,2)$ for $e=0$ and $e\geq 2$. In Section~\ref{sec:Hyper-bitangent Curves on F1}, we examine $\operatorname{Hyp}(\mathbb{F}_{1},B,2)$ for the non-minimal Hirzebruch surface $\mathbb{F}_{1}\cong \operatorname{Bl}_{pt} \mathbb{P}^{2}$.

\subsection*{Notations} We work over $\mathbb{C}$ unless otherwise stated. Our notations will be consistent with \cite{caporaso2024hypertangencyplanecurvesalgebraic}.

\begin{itemize}
    \item Throughout the paper, $S$ is an irreducible smooth surface, curves are assumed to be reduced, and curves on a surface $S$ are assumed to be closed in $S$.
    \item For a divisor $D$ on $S$, we denote by $\mathcal{L}(D)$ the associated divisor class of it in $\operatorname{Pic}(S)$.
    \item Let $B$ and $C$ be two curves on $S$, and $q \in B \cap C$ be a point. We say $B$ and $C$ intersect transversally at $q$ if their locally irreducible branches at $q$ share no common tangent directions; otherwise, we say they are tangent at $q$.
    \item For a $3C$-curve $B$ on a smooth projective surface $S$, we denote by $N$ the set of intersection points of its components. Write $B=B_{1}\cup B_{2}\cup B_{3}$ where the $B_{i}$'s are the irreducible components of $B$, then $N=\bigcup_{i < j} B_{i}\cap B_{j}$. And by $p_{i,j}$, we mean an intersection point of $B_{i}$ and $B_{j}$.
    \item We denote by $\operatorname{mult}_{p}(D)$ the multiplicity of a point $p$ on a curve $D$.
    \item For an integral curve $D$, we denote by $\nu_{D}\colon D^{\nu}\to D$ the normalization map of $D$.
    \item We say a curve $D$ is a rational curve if its normalization has arithmetic genus $p_{a}(D^{\nu}) = 0$.
    \item For a point $p$ on an integral curve $D$, its $\delta$-invariant is defined as
    $
    \delta_{D}(p) \coloneqq p_{a}(D) - p_{a}(D_{p}^{\nu})
    $,
    where $D_{p}^{\nu}$ is the partial normalization of $D$ at $p$. It is well-known that 
    $
    \delta_{D}(p) = \sum_{p^{\prime}} \frac{m_{p^{\prime}}(m_{p^{\prime}}-1)}{2}
    $,
    where $p^{\prime}$ runs over all the infinitely near points lying over $p$, including $p$ itself, and $m_{p^{\prime}}$ is the multiplicity of $p^{\prime}$. 
    \item Let $B$ be a a curve on a smooth projective surface $S$. We say an integral curve $D$ on $S$ is \textit{hypertangent} to $B$ if
    $\#\nu_{D}^{-1}(D\cap B) = 1$;
    we say $D$ is \textit{hyper-bitangent} to $B$ if 
    $\#\nu_{D}^{-1}(D\cap B) = 2$.
    \item Let $B$ be a curve on a smooth projective surface $S$, and $\mathcal{C}\in \operatorname{Pic}(S)$ be a divisor class. We define
    $$
        \operatorname{Hyp}_{\mathcal{C}}(B,2) = 
        \left\{
            D \in |\mathcal{C}| 
            \ \big\vert \
            D \text{ is an integral curve such that }\#\nu_{D}^{-1}(D \cap B) \leq 2
        \right\}.
    $$
    Clearly, $\operatorname{Hyp}(B,2) = \bigcup_{\mathcal{C}\in \operatorname{Pic(S)}} \operatorname{Hyp}_{\mathcal{C}}(B,2)$.
    
    For a point $p\in B$, we define
    $$
    \operatorname{Hyp}_{\mathcal{C}}(B,p) 
    =
    \left\{
    D \in |\mathcal{C}| 
    \ \big\vert \
    D \text{ is an integral curve which is } \text{hypertangent to $B$ at $p$}
    \right\},
    $$
    and
    $\operatorname{Hyp}_{\mathcal{C}}^{m}(B,p) 
    =
    \left\{
    D \in \operatorname{Hyp}_{\mathcal{C}}(B,p) 
    \ \big\vert \
    \operatorname{mult}_{p}(D) = m
    \right\}$.
\end{itemize}

\subsection*{Acknowledgements}
The author would like to express gratitude to his supervisor, Lucia Caporaso, for suggesting this problem, and for numerous discussions and constant encouragement throughout the development of this article. The author thanks Amos Turchet for helpful conversations, especially for pointing out the references \cite{CZ08} and \cite{CZ13}. Additionally, the author thanks Davide Gori, Riccardo Moschetti, Lidia Stoppino, and Roberto Vacca for helpful discussions. Lastly, the author would also like to thank the referees for their helpful comments and suggestions.

\section{Preliminaries}
\label{sec:prelim}

\subsection{Unibranch Points of Algebraic Curves}\label{subsec:unibranch points}\hfill

All the results in this subsection hold over an algebraically closed field $\Bbbk$ of arbitrary characteristic.

\subsubsection{Analytic Invariants of Algebraic Curves}\hfill

In this sub-subsection, we firstly collect some facts concerning analytic invariants of algebraic curves. Here by "analytic" we mean the invariant is determined by the completion of the local ring involved. All the results in this sub-subsection are known.

\begin{lemma}\label{lem: multiplicity is analytic invariant}
    Let $\varphi \colon C \to D$ be a morphism of curves. If $\varphi$ induces an isomorphism of completions of local rings on $q\in C$ and $\varphi(q)\in D$, that is, $\varphi^{*} \colon \hat{\mathcal{O}}_{\varphi(q),D}\to \hat{\mathcal{O}}_{q,C}$ is an isomorphism, then
    $\operatorname{mult}_{q}(C)=\operatorname{mult}_{\varphi(q)}(D)$. 
\end{lemma}

\begin{lemma}\label{lem: delta-invariant is locally analytic}
    The $\delta-$invariant of a point on a curve is an analytic invariant. Let $\varphi\colon C\to D$ be a morphism of curves. If $\varphi$ induces an isomorphism of completions of local rings on $q\in C$ and $\varphi(q)\in D$, that is, $\varphi^{*}\colon \hat{\mathcal{O}}_{\varphi(q),D}\to \hat{\mathcal{O}}_{q,C}$ is an isomorphism, then
    $\delta_{C}(q)=\delta_{D}(\varphi(q))$. 
\end{lemma}
\begin{proof}
    See \cite[\href{https://stacks.math.columbia.edu/tag/0C3Q}{Tag 0C3Q}]{stacks-project} and \cite[\href{https://stacks.math.columbia.edu/tag/0C1R}{Tag 0C1R}]{stacks-project}.
\end{proof}

\begin{lemma}\label{lem: number of branches of a curve is locally analytic}
    Let $q$ be a point on a curve $C$; then the number of branches of $C$ at $q$ is an analytic invariant.
\end{lemma}
\begin{proof}
    See \cite[\href{https://stacks.math.columbia.edu/tag/0C2D}{Tag 0C2D}]{stacks-project} and \cite[\href{https://stacks.math.columbia.edu/tag/0C3Z}{Tag 0C3Z}]{stacks-project}.
\end{proof}

\begin{lemma}\label{lem: Intersection multiplicity of two curves on a surface is locally analytic}
    Let $C, D$ be two curves on a surface $S$, and let $q$ be a point in $C \cap D$. Then
    $(C \cdot D)_q = \dim_{\Bbbk} \widehat{\mathcal{O}}_{q,S}/(f_C, f_D)$.
    where $f_C$ (respectively $f_D$) is a local equation defining $C$ (respectively $D$).
\end{lemma}

\subsubsection{Properties of Unibranch Points}\hfill

In this sub-subsection we study unibranch points of algebraic curves.

We firstly introduce the notion of \textit{multiplicity sequence} due to Flenner--Zaidenberg \cite{FlennerZaidenberg1996}:
\begin{definition}
    Let $q$ be a singular unibranch point on a curve $C$, which is contained in a surface $S$. Let 
    $$
    S = V_{0} 
    \xleftarrow{\sigma_{1}} V_{1}
    \xleftarrow{\sigma_{2}}
    \cdots
    \xleftarrow{\sigma_{n}} V_{n}
    $$
    be a minimal resolution of $q$, $C^{i}$ be the strict transform of $C$ in $V_{i}$ and $q_{i}\in C^{i}$ be the point lying over $q$; we also denote $C^{0}=C$ and $q_0 = q$.
    
    Denote by $m_{i}$ the multiplicity of $q_{i}\in C^{i}$, then the multiplicity sequence of $q$ is defined to be the following sequence of positive integers:
    $$
    \underline{m}_{q} = (m_{0},m_{1},\cdots,m_{n}).
    $$
    
    We denote $m \coloneqq \operatorname{mult}_{q}(C)$ and by $l(q)$ the smallest integer such that $m_{l(q)} < m$. With these notations, we have $m = m_{0} = \cdots = m_{l(q)-1}$ and $m_n = 1$.
\end{definition}

Then we have the following lemma by Flenner--Zaidenberg \cite[Lem.~1.4]{FlennerZaidenberg1996} describing the intersection multiplicity of $C$ with another curve which is smooth at $q$:
\begin{lemma}\label{lem: Flenner-Zaidenberg intersection multiplicity at a singular unibranch point}
    Let $B, C$ be two curves in $\mathbb{A}^{2}$, and let $q \in B\cap C$ be a unibranch point for both $B$ and $C$. If $q$ is a smooth point of $B$ and $\operatorname{mult}_{q}(C) = m \geq 2$, then either
    $
    (B\cdot C)_{q} = k m
    $
    for some integer $k$ with $1 \leq k \leq l(q)$,
    or
    $
    (B\cdot C)_{q} = l(q) m + m_{l(q)}
    $.
\end{lemma}

\begin{remark}\label{rmk: Flenner-Zaidenberg intersection multiplicity at a singular unibranch point}
    By our previous study of analytic invariants, we observe that Lem.~\ref{lem: Flenner-Zaidenberg intersection multiplicity at a singular unibranch point} actually holds over any smooth surface $S$, as we will see in the proof of Thm.~\ref{thm: bounddeltainvariant}. 
\end{remark}

\begin{remark}\label{rmk: more on Flenner-Zaidenberg intersection multiplicity at a singular unibranch point}
    In \cite[Lem.~1.4]{FlennerZaidenberg1996}, Flenner and Zaidenberg established a stronger result: on $\mathbb{A}^2$, given a singular unibranch point $q$ of a curve $C$, for any possible intersection multiplicity described in Lem.~\ref{lem: Flenner-Zaidenberg intersection multiplicity at a singular unibranch point}, there exists a curve $B$ containing $q$ as a smooth point such that $(B \cdot C)_{q}$ attains the given intersection multiplicity.

    In fact, this works over any surface $S$: let $q$ be a singular unibranch point of a curve $C \subset S$. Since $S$ is smooth, by \cite[Prop.~4.9]{milneLEC}, there exists an open neighborhood $U$ of $q$ and a regular morphism $\varphi \colon U \to \mathbb{A}^{2}$ that is étale at $q$. From Rem.~\ref{rmk: Flenner-Zaidenberg intersection multiplicity at a singular unibranch point}, we know $\varphi(q)$ is a singular unibranch point of $\varphi(C)$ with 
    $\underline{m}_{\varphi(q)} = \underline{m}_q$. By \cite[Lem.~1.4]{FlennerZaidenberg1996}, we can find a curve $B \subset \mathbb{A}^2$ such that $B$ contains $\varphi(q)$ as a smooth point and intersects $\varphi(C)$ at $\varphi(q)$ with the prescribed intersection multiplicity. Then $\varphi^{-1}(B)$ contains $q$ as a smooth point and intersects $C$ at $q$ with the same intersection multiplicity. 
\end{remark}

\begin{corollary}\label{cor: Intersection Multiplicity is Invariant for (m,n)-point with n<2m}
Let $q \in C \subset \mathbb{A}^{2}$ be a singular unibranch point with $\operatorname{mult}_{q}(C)=m \geq 2$, and let $L$ be the tangent line to $C$ at $q$. If $(C\cdot L)_{q} = n < 2m$, then for any curve $B \subset \mathbb{A}^{2}$ containing $q$ as a smooth point such that $B$ is tangent to $C$ at $q$, we have
$
(B\cdot C)_{q} = n
$,
\end{corollary}

\begin{proof}
By assumption we have $m < n < 2m$. From Lem.~\ref{lem: Flenner-Zaidenberg intersection multiplicity at a singular unibranch point} we see that in the multiplicity sequence $\underline{m}_{q}$, we have $m_{0} = m$ and $m_{1} = n - m < m$. Hence, we conclude that $l(q) = 1$, and the result follows immediately.
\end{proof}

\begin{remark}
    Let $C \subset \mathbb{P}^{2}$ be the curve defined by the equation $z^{n-m}y^{m} - x^{n} = 0$ with $m < n < 2m$ and $m\geq 2$. We see that $q = [0:0:1]$ is a singular unibranch point of $C$ with $\operatorname{mult}_{q}(C) = m$. Let $L$ be the tangent line to $C$ at $q$, then $L$ is defined by $y=0$ and $(C\cdot L)_{q} = n$. If $B$ is a curve that contains $q$ as a smooth point and is hypertangent to $C$ at $q$, then, by Cor.~\ref{cor: Intersection Multiplicity is Invariant for (m,n)-point with n<2m}, we have  
    $
    n = (B \cdot C)_{q} = (B \cdot C) = n \deg B
    $.
    This gives $\deg B = 1$, so $B$ must be the tangent line to $C$ at $q$, meaning $B = L$.
\end{remark}

The following theorem will be one of the key tools we will use later:

\begin{theorem}\label{thm: bounddeltainvariant}
Let $B, C$ be two curves on a surface $S$, and $q \in B \cap C$ be a unibranch point for both $B$ and $C$. If $q$ is a smooth point of $B$ and $\operatorname{mult}_{q}(C) = m$, then  
\begin{equation}
    \delta_{C}(q) \geq \frac{(m-1) \left( (B\cdot C)_{q} - 1 \right) + \operatorname{gcd} \left( (B\cdot C)_{q}, m \right) - 1}{2} 
    \geq \frac{(m-1) \left( (B\cdot C)_{q} - 1 \right)}{2},
    \nonumber
\end{equation}
and the two equalities hold if $\operatorname{gcd} \left( (B\cdot C)_{q}, m \right) = 1$.
\end{theorem}

\begin{proof}
    Since $S$ is smooth, by \cite[Prop.~4.9]{milneLEC}, there exists an open neighborhood $U$ of $q$ and a regular morphism $\varphi \colon U \to \mathbb{A}^{2}$ that is étale at $q$.
    
    Let $f$ and $g$ be the defining local equations for $\varphi(B)$ and $\varphi(C)$, respectively. Then, $\varphi^{*}f$ and $\varphi^{*}g$ are the local equations defining $B$ and $C$, respectively. Since $\varphi$ is étale at $q$, it induces an isomorphism
    $
    \varphi^{*} \colon \hat{\mathcal{O}}_{\varphi(q), \mathbb{A}^{2}} \to \hat{\mathcal{O}}_{q, S}
    $.
    Combining Lem.~\ref{lem: Intersection multiplicity of two curves on a surface is locally analytic}, we obtain
    $$
    \mathcal{O}_{\varphi(q), \mathbb{A}^{2}}/(f,g) \cong \hat{\mathcal{O}}_{\varphi(q), \mathbb{A}^{2}}/(f,g) \cong \hat{\mathcal{O}}_{q, S}/(\varphi^{*}f, \varphi^{*}g) \cong \mathcal{O}_{q, S}/(\varphi^{*}f, \varphi^{*}g).
    $$
    Hence, we have
    $
    (B \cdot C)_{q} = (\varphi(B) \cdot \varphi(C))_{\varphi(q)}
    $.
    Moreover, by Lem.~\ref{lem: multiplicity is analytic invariant} and Lem.~\ref{lem: number of branches of a curve is locally analytic}, we know that $\varphi(q)$ is a smooth point of $\varphi(B)$ and a unibranch $m$-fold point of $\varphi(C)$. By Lem.~\ref{lem: delta-invariant is locally analytic}, we know that $\delta_C(q) = \delta_{\varphi(C)}(\varphi(q))$. Thus, $\varphi(B)$, $\varphi(C)$, and $\varphi(q)$ satisfy all the assumptions in the statement of the theorem. Therefore, we may assume $S = \mathbb{A}^2$.

    Now, in the affine plane, if both $B$ and $C$ are lines, then $m = 1$ and $\delta_C(q) = 0$, so the statement is trivially true.

    If $m = 1$, we have $\delta_C(q) = 0$, the statement is also trivially true.
    
    If $B$ intersects $C$ transversally at $q$, then $(B \cdot C)_{q} = m$, and we have $\delta_C(q) \geq \frac{m(m-1)}{2} \geq \frac{(m-1)^2}{2}$, as desired.

    From now on, we assume $m \geq 2$ and $B$ is tangent to $C$ at $q$, so $(B \cdot C)_{q} > m$. Assume $(B \cdot C)_{q} = km + r$, where $k, r \in \mathbb{N}$, $k \geq 1$, and $m > r \geq 0$.

    If $r = 0$, then
    $
    \delta_C(q) \geq \left(l(q) m (m-1)\right)/2 \geq \left((B \cdot C)_{q} \cdot (m-1)\right)/2
    $
    as desired.

    If $r > 0$, then by Lem.~\ref{lem: Flenner-Zaidenberg intersection multiplicity at a singular unibranch point}, we know that $k = l(q)$.
    
    If $k = 1$ and $r > 0$, let $L_q$ be the tangent line to $C$ at $q$. By Cor.~\ref{cor: Intersection Multiplicity is Invariant for (m,n)-point with n<2m}, we know that $(B \cdot C)_{q} = (L_q \cdot C)_{q} = m + r$, so the result follows from \cite[Prop.~5.1.1]{caporaso2024tropicalcurvesunibranchpoints}.

    If $k \geq 2$ and $r > 0$, since $l(q) = k \geq 2$, we know that $\operatorname{mult}_{q_1}(C^1) = m$ and $\underline{m}_{q_1} = (m_1, \dots, m_n)$. Then, we can proceed by induction on $k$, we obtain
    \begin{equation}
        \begin{split}
            \delta_C(q) &= \delta_{C^1}(q_1) + \frac{m(m-1)}{2} \\
            &\geq \frac{(m-1)\left( (k-1)m + r - 1 \right) + \operatorname{gcd}\left( (k-1)m + r, m \right) - 1}{2} + \frac{m(m-1)}{2} \\
            &= \frac{(m-1)\left( km + r - 1 \right) + \operatorname{gcd}(r, m) - 1}{2} \\
            &= \frac{(m-1)\left( (B \cdot C)_{q} - 1 \right) + \operatorname{gcd}\left( (B \cdot C)_{q}, m \right) - 1}{2}
        \end{split}
    \end{equation}
    as desired. This completes the proof.
\end{proof}

\begin{remark}
    Thm.~\ref{thm: bounddeltainvariant} is the generalized local version of the lower bound of $\delta$-invariant in \cite[Thm.~2.2.1]{caporaso2024hypertangencyplanecurvesalgebraic}.
\end{remark}

Another key tool we will also need is the following \textbf{Strong Triangle Inequality} due to Garc\'ia Barroso--Płoski \cite[Thm.~2.8, Cor.~2.9]{BP15}:
\begin{theorem}\label{thm: Strong Triangle Inequality}
    Let $B,C,D$ be three curves on a surface $S$. If they intersect at a point $q$ which is unibranch on all of them, then the smallest two among
    $$
    \frac{(B\cdot C)_{q}}{\operatorname{mult}_{q}(B)\operatorname{mult}_{q}(C)},
    \quad
    \frac{(B\cdot D)_{q}}{\operatorname{mult}_{q}(B)\operatorname{mult}_{q}(D)},
    \quad
    \frac{(C\cdot D)_{q}}{\operatorname{mult}_{q}(C)\operatorname{mult}_{q}(D)}
    $$
    are equal.
\end{theorem}

\begin{remark}
    We note that the strong triangle inequality in Thm.~\ref{thm: Strong Triangle Inequality} over $\mathbb{C}$ was first proved by P{\l}oski in \cite{Ploski1985}. 
\end{remark}

\subsection{Basics on Hirzebruch Surfaces}\label{subsec:Hirzebruch surfaces}\hfill

Recall that the Hirzebruch surface $\mathbb{F}_{e}$ is $\mathbb{P}_{\mathbb{P}^{1}}(\mathcal{O} \oplus \mathcal{O}(e))$ for $e \geq 0$. We collect some basic facts about Hirzebruch surfaces here for convenience, with notations consistent with \cite[Ch.~\Romannum{5}, Section 2]{Hartshorne1977}.

Let $\rho \colon \mathbb{F}_{e} \to \mathbb{P}^1$ denote the ruling. We denote by $f$ the class of a fiber of $\rho$. There are two disjoint sections: $C_0$, the section with self-intersection $-e$, and $C_1$, the section with self-intersection $e$. In particular, when $e=0$, we have $C_0 \sim C_1$, which corresponds to the fibers of the second ruling on $\mathbb{F}_0 \cong \mathbb{P}^1 \times \mathbb{P}^1$.

\begin{theorem}
    We have $\operatorname{Pic}(\mathbb{F}_{e}) = \mathbb{Z} C_{1} \oplus \mathbb{Z} f$, and 
    \begin{itemize}
        \item $C_{1}^{2} = e$, $C_{0}^{2} = -e$, $f^{2} = 0$;
        \item $(C_{0} \cdot f) = 1 = (C_{1} \cdot f)$, $(C_{0} \cdot C_{1}) = 0$;
        \item $C_{0} \sim C_{1} - e f$, $K_{\mathbb{F}_{e}} \sim -2 C_{1} + (e - 2) f$.
    \end{itemize}
\end{theorem}
\begin{proof}
    See \cite[Ch.~\Romannum{5}, Lem.~2.10, p.~373 and Thm.~2.17, p.~379]{Hartshorne1977}.
\end{proof}

\begin{corollary}\label{cor : Ample divisor on Hirzebruch surface}
    Let $\alpha, \beta \in \mathbb{Z}$. Then the following are equivalent:
    \begin{enumerate}
        \item $\alpha C_{1} + \beta f$ is very ample;
        \item $\alpha C_{1} + \beta f$ is ample;
        \item $\alpha > 0$ and $\beta > 0$.
    \end{enumerate}
\end{corollary}

\begin{corollary}\label{cor: Irreducible Curves on Hirzebruch Surfaces}
    Let $\alpha, \beta \in \mathbb{Z}$. Then the following are equivalent:
    \begin{enumerate}
        \item $\lvert \alpha C_{1} + \beta f \rvert$ contains a smooth irreducible curve;
        \item $\lvert \alpha C_{1} + \beta f \rvert$ contains an integral curve;
        \item $\alpha = 0$, $\beta = 1$, or $\alpha = 1$, $\beta = -e$, or $\alpha > 0$, $\beta \geq 0$.
    \end{enumerate}
\end{corollary}

\begin{proof}(of the two corollaries)
    See \cite[Ch.~\Romannum{5}, Cor.~2.18, p.~380]{Hartshorne1977}.
\end{proof}

\begin{remark}\label{rmk: Effective divisors on Hirzebruch Surfaces}
    From Cor.~\ref{cor: Irreducible Curves on Hirzebruch Surfaces}, we easily see that a divisor class $u C_{1} + v f = u C_{0} + (u e + v) f \in \operatorname{Pic}(\mathbb{F}_{e})$ is effective if and only if $u \geq 0$ and $u e + v \geq 0$, because an effective divisor is a formal sum of integral divisors with non-negative coefficients.
\end{remark}

\begin{lemma}\label{lem: Global section of Effective divisor on Hirzebruch Surfaces}
    For an effective divisor class $u C_{1} + v f \in \operatorname{Pic}(\mathbb{F}_{e})$ ($u \geq 0$, $u e + v \geq 0$), 
    if $e \geq 1$, we have 
    $$
    h^{0}(\mathbb{F}_{e}, u C_{1} + v f)
    = h^0 \left( \mathbb{F}_{e}, u C_{0} + (u e + v) f 
    \right)
    = \sum_{i = 0}^{\min \{u, \lfloor \frac{u e + v}{e} \rfloor\}} (u e + v - i e + 1);
    $$
    
    if $e = 0$, we have
    $
    h^{0}(\mathbb{F}_{e}, u C_{1} + v f) = (u + 1)(v + 1)
    $.
\end{lemma}
\begin{proof}
    See \cite{MO403548}.
\end{proof}

\begin{remark}\label{rmk: C0 unique negative self intersection curve on Fe}
When $e \geq 1$, we have $C_{0} \sim C_{1} - e f$, and by Lem.~\ref{lem: Global section of Effective divisor on Hirzebruch Surfaces}, we obtain $h^{0}(\mathbb{F}_{e}, C_{0}) = 1$. Therefore, there exists a unique integral curve in $\lvert C_{0} \rvert$, which is the unique smooth rational curve on $\mathbb{F}_{e}, e\geq 1$, with negative self-intersection number. We also denote this curve as $C_{0}$.
\end{remark}

The following lemma, which follows directly from the adjunction formula, is well known:
\begin{lemma}
    On $\mathbb{F}_{e}$, an integral curve $D \in \lvert \alpha C_{1} + \beta f \rvert$ has arithmetic genus
    $$
    p_{a}(D)
    = \frac{1}{2} e \alpha (\alpha - 1) + (\alpha - 1)(\beta - 1)
    = \frac{1}{2} (\alpha - 1)(e \alpha + 2 \beta - 2).
    $$
\end{lemma}

\begin{remark}\label{rmk: F1 is blow-up a point on P2}
    Let us determine when $p_{a}(D) = 0$. There are two cases: $\alpha = 1$ and $e \alpha + 2 \beta - 2 = 0$.

    The first case, $\alpha = 1$, tells us that $\lvert C_{1} + \beta f \rvert$ with $\beta \geq 0$ is a linear system whose irreducible members are smooth rational curves.
    
    For the second case, we have $2 = e \alpha + 2 \beta$. If $\beta > 0$, then $2 = e \alpha + 2 \beta \geq 2 \beta \geq 2$, hence the equalities here hold; we must have $\alpha = 0$ and $\beta = 1$, in this case $D \in \lvert f \rvert$ is a fiber. If $\beta = 0$, then $e \alpha = 2$ splits into two subcases: either $e = 1$, $\alpha = 2$ or $e = 2$, $\alpha = 1$. Since the latter case is already included in the $\alpha = 1$ case, the only remaining case is $e = 1$, $\alpha = 2$.
    
    When $e = 1$, $\mathbb{F}_{1} \cong \operatorname{Bl}_{p} \mathbb{P}^{2}$, where $p$ is a point in $\mathbb{P}^{2}$. Denote the blow-up morphism by $\pi \colon \operatorname{Bl}_{p} \mathbb{P}^{2} \to \mathbb{P}^{2}$ and the exceptional divisor by $E$. We have the following correspondences:
    $$
    C_{1} = \pi^{*} \mathcal{O}(1),
    \quad
    C_{0} = E,
    \quad
    f = \pi^{*} \mathcal{O}(1) - E.
    $$
    
    We see that the only special case (given by $e \alpha + 2 \beta - 2 = 0$) for smooth rational curves on $\mathbb{F}_{e}$ is $\pi^{*} \mathcal{O}(2) \in \operatorname{Pic}(\mathbb{F}_{1})$.
\end{remark}

\begin{lemma}\label{lem: Upper bound of unibranch multiplicity on Hirzebruch Surfaces}
    For a unibranch point $q$ of an integral curve $C \subset \mathbb{F}_{e}$, assume $C \in \lvert \alpha C_{1} + \beta f \rvert$ and $C$ is not in $\lvert C_{0} \rvert$ or $\lvert f \rvert$, and let $m = \operatorname{mult}_{q}(C)$.
    
    If $e = 0$, then $\alpha, \beta > 0$ and $m \leq \min\{\alpha, \beta\}$.

    If $e \geq 1$, then $\alpha > 0$, $\beta \geq 0$, and $m\leq \alpha$.
\end{lemma}
\begin{proof}
    If $e = 0$, then $C_{1} \sim C_{0}$, and $\lvert C_0 \rvert = \lvert C_{1} \rvert$ and $\lvert f \rvert$ are two rulings on $\mathbb{F}_{0}$. Hence, $C \not\in \lvert C_0 \rvert = \lvert C_1 \rvert$ and $C \not\in \lvert f\lvert$ imply that $\alpha, \beta > 0$, and there exist unique $C_{1,q} \in \lvert C_{1} \rvert$ and $f_{q} \in \lvert f \rvert$ that contain $q$. Thus, we obtain 
    $$
    m
    \leq \min \{(f_{q} \cdot C)_{q}, (C_{1,q} \cdot C)_{q}\}
    \leq \min \{(f_{q} \cdot C), (C_{1,q} \cdot C)\}
    = \min \{\alpha, \beta\}
    $$
    as desired.

    If $e \geq 1$, then by Cor.~\ref{cor: Irreducible Curves on Hirzebruch Surfaces}, we know the assumption $C \not\in \lvert C_0\rvert \text{ or } \lvert f \rvert$ implies that $\alpha > 0$ and $\beta \geq 0$. In this case, we note that either there exists an integral curve $L \in \lvert C_{1} \rvert$ such that $q \in L$ and $L \neq C$, or $q \in C_{0}$.

    Let $f_{q} \in \lvert f \rvert$ be the unique fiber containing $q$, then
    $$
    m
    \leq (f_{q} \cdot C)_{q}
    \leq (f_{q} \cdot C)
    = \left( f \cdot (\alpha C_{1} + \beta f) \right)
    = \alpha.
    $$

    If $q \in L \in \lvert C_{1} \rvert$, then
    $$
    m
    \leq (L \cdot C)_{q}
    \leq (L \cdot C)
    = \left( C_{1} \cdot (\alpha C_{1} + \beta f) \right)
    = e \alpha + \beta.
    $$
    In this case, we get $m \leq \alpha$.

    If $q \in C_{0}$, then
    $$
    m
    \leq (C_{0} \cdot C)_{q}
    \leq (C_{0} \cdot C)
    = \left( C_{0} \cdot (\alpha C_{1} + \beta f) \right)
    = \beta.
    $$
    In this case, we get $m \leq \min\{\alpha, \beta\} \leq \alpha$.
\end{proof}

\begin{lemma}\label{lem: Special case of bound of unibranch multiplicity on F1}
    For a point $q$ of an integral curve $C \subset \mathbb{F}_{1}$ where $C \in \lvert d C_{1} \rvert$ for some $d \geq 2$, we have 
    $
    \operatorname{mult}_{q}(C) < d
    $.
\end{lemma}
\begin{proof}
    This follows from the fact that $d C_{1} = \pi^{*} \mathcal{O}(d)$, as mentioned in Rem.~\ref{rmk: F1 is blow-up a point on P2}, and an integral plane curve of degree $d\geq 2$ cannot have a point of multiplicity $\geq d$.
\end{proof}

\begin{definition}
    Let $X$ be an irreducible projective variety of dimension $m$, and let $\mathcal{L} \in \operatorname{Pic}(X)$. The volume of $\mathcal{L}$ is defined as 
    $$
    \operatorname{vol}(\mathcal{L}) = \limsup_{n \to \infty} \frac{h^{0}(X, \mathcal{L}^{\otimes n})}{n^{m} / m!}.
    $$
    Let $D$ be a Cartier divisor on $X$, we denote by $\mathcal{L}(D)$ the associated divisor class of it in $\operatorname{Pic}(X)$. The volume of $D$ is defined as $\operatorname{vol}(\mathcal{L}(D))$.

    It is well-known that $\mathcal{L}$ is big if and only if $\operatorname{vol}(\mathcal{L}) > 0$; see, for instance, \cite[Ch.~2, Section~2.2.C, p.~148]{Lazarsfeld2004}.
\end{definition}

\begin{proposition}\label{prop: Criterion for bigness on Hirzebruch Surfaces}
    A divisor class $u C_{1} + v f \in \operatorname{Pic}(\mathbb{F}_{e})$ is big if and only if $u > 0$ and $u e + v > 0$.
\end{proposition}
\begin{proof}
    First, we note that a big divisor class on $\mathbb{F}_{e}$ must be effective. Indeed, if $u C_{1} + v f$ is a big divisor class on $\mathbb{F}_{e}$, then there exist $m \in \mathbb{N}$, an ample divisor class $A$, and an effective divisor class $F$ such that $m\left(u C_{1} + v f\right) \sim m\left(u C_{0} + (u e + v) f\right) \sim A + F$. By Cor.~\ref{cor : Ample divisor on Hirzebruch surface} and Rem.~\ref{rmk: Effective divisors on Hirzebruch Surfaces}, we obtain $m u > 0$ and $m(u e + v) > 0$, so $u > 0$ and $u e + v > 0$. Thus, $u C_{1} + v f$ is effective.

    For an effective divisor class $u C_{1} + v f$ ($u \geq 0$, $u e + v \geq 0$), we use Lem.~\ref{lem: Global section of Effective divisor on Hirzebruch Surfaces} to compute its volume $\operatorname{vol}(u C_{1} + v f)$.

    If $e = 0$, we have 
    $
    h^{0}(\mathbb{F}_{0}, n(u C_{1} + v f))
    = (n u + 1)(n v + 1)
    = u v n^{2} + (u + v) n + 1
    $.
    Thus
    $$
        \operatorname{vol}(u C_{1} + v f)
        = \limsup_{n \to \infty} \frac{h^{0}(\mathbb{F}_{0}, n(u C_{1} + v f))}{n^{2} / 2}
        = \limsup_{n \to \infty} \frac{2 (u v n^{2} + (u + v) n + 1)}{n^{2}}
        = 2 u v.
    $$
    Hence, $u C_{1} + v f$ is big if and only if $u > 0$ and $v > 0$.

    If $e \geq 1$, we consider the following subcases:
    \begin{enumerate}
        \item If $u = 0$, then $h^{0}(\mathbb{F}_{e}, n(u C_{1} + v f)) = h^{0}(\mathbb{F}_{e}, n v f) = n v + 1$. In this subcase, $u C_{1} + v f = v f$ is not big.
        
        \item If $u > 0$ and $v \geq 0$, then
        \begin{equation}
            \begin{split}
                h^{0}(\mathbb{F}_{e}, n(u C_{1} + v f)) 
                & = \sum_{i = 0}^{n u} (n u e + n v - i e + 1) \\
                & = \left( (u e + v) n + 1 \right) (u n + 1) - e \sum_{i = 0}^{u n} i \\
                & = u (u e + v) n^{2} + (u e + v + u) n + 1 - \frac{e u^{2} n^{2} + e u n}{2} \\
                & = \left( u v + \frac{e u^{2}}{2} \right) n^{2} + \left( u + v + \frac{e u}{2} \right) n + 1.
            \end{split}
            \nonumber
        \end{equation}
        Thus
        $$
            \operatorname{vol}(u C_{1} + v f)
            = \limsup_{n \to \infty} \frac{ \left( u v + \frac{e u^{2}}{2} \right) n^{2} + \left( u + v + \frac{e u}{2} \right) n + 1 }{n^{2} / 2}
            = 2 u v + e u^{2} > 0.
        $$
        Hence, $u C_{1} + v f$ is big when $u > 0$ and $v \geq 0$.
        
        \item If $u > 0$ and $v < 0$, we write $\lfloor \frac{n u e + n v}{e} \rfloor = \frac{n u e + n v}{e} - \epsilon$, where $\epsilon$ is a rational number determined by $n$, $v$, $e$, and $0 \leq \epsilon < 1$. Then
        \begin{equation}
            \begin{split}
                h^{0}(\mathbb{F}_{e}, n(u C_{1} + v f)) 
                & = \sum_{i = 0}^{\frac{n u e + n v}{e} - \epsilon} (n u e + n v - i e + 1) \\
                & = \left( (u e + v) n + 1 \right) \left( \frac{(u e + v) n}{e} - \epsilon \right) - e \sum_{i = 0}^{\frac{n u e + n v}{e} - \epsilon} i \\
                & = \frac{(u e + v)^{2} n^{2}}{e} + \frac{(u e + v)(1 - e \epsilon) n}{e} + \epsilon \\
                & \quad - \frac{e}{2} \left( \frac{(u e + v)^{2} n^{2}}{e^{2}} + \frac{(1 - 2 \epsilon)(u e + v) n}{e} \right) - \frac{e}{2} (\epsilon^{2} - \epsilon) \\
                & = \frac{(u e + v)^{2}}{2 e} n^{2} + \frac{(u e + v)(2 - e)}{2 e} n + \epsilon - \frac{e}{2} (\epsilon^{2} - \epsilon).
            \end{split}
            \nonumber
        \end{equation}
        Thus
        $$
            \operatorname{vol}(u C_{1} + v f)
            = \limsup_{n \to \infty} \frac{ \frac{(u e + v)^{2}}{2 e} n^{2} + \frac{(u e + v)(2 - e)}{2 e} n + \epsilon - \frac{e}{2} (\epsilon^{2} - \epsilon) }{n^{2} / 2}
            = \frac{(u e + v)^{2}}{e} \geq 0.
        $$
        Since $u e + v \geq 0$, we see that $u C_{1} + v f$ is big if and only if $u e + v > 0$ in this subcase.
    \end{enumerate}
    In summary, an effective divisor class $u C_{1} + v f$ on $\mathbb{F}_{e}$ is big if and only if $u > 0$ and $u e + v > 0$ when $e \geq 1$. Combining this with the $e = 0$ case, the proof is complete.
\end{proof}

\section{Hyper-bitangent Curves on Minimal Hirzebruch Surfaces}
\label{sec:Hyper-bitangent Curves on Minimal Hirzebruch Surfaces}
In this section, we investigate $\operatorname{Hyp}(B,2)$ for a minimal Hirzebruch surface \( S \), that is, $S = \mathbb{F}_{e}$ for $e = 0$ and $e \geq 2$. 

On $\mathbb{F}_{0} \cong \mathbb{P}^{1} \times \mathbb{P}^{1}$, we denote the divisor $u C_{1} + v f$ as $(u, v)$ by tradition. Notice that $C_{0}=C_{1}$ in this case. With the polarization $\left( \mathbb{P}^{1} \times \mathbb{P}^{1}, (1,1) \right)$ fixed, $\lvert C_1 \rvert = \lvert (1,0) \rvert$ and $\lvert f \rvert = \lvert (0,1) \rvert$ give two families of lines. By Prop.~\ref{prop: Criterion for bigness on Hirzebruch Surfaces}, we see that for a divisor $B$ on $\mathbb{P}^{1} \times \mathbb{P}^{1}$, $ K_{ \mathbb{P}^{1} \times \mathbb{P}^{1} } + \mathcal{L}(B)$ is big if and only if $\mathcal{L}(B)\geq (3,3)$.

Then, we have the following:

\begin{proposition}\label{prop: F0 Hyp(B,2)=E(B)}
    Consider $\mathbb{F}_{0} \cong \mathbb{P}^{1}\times\mathbb{P}^{1}$. Let $B=B_{1}\cup B_{2}\cup B_{3} \in \lvert \alpha C_{1} + \beta f \rvert$ be a $3C$-curve such that $K_{\mathbb{P}^{1} \times \mathbb{P}^{1}} + \mathcal{L}(B)$ is big. Assume $B_{i} \in \lvert(\alpha_{i}, \beta_{i})\rvert$ with $(\alpha_i , \beta_i) \geq (1, 1)$ for all $i = 1, 2, 3$, and $\alpha_1 + \beta_1 \leq \alpha_2 + \beta_2 \leq \alpha_3 + \beta_3$. Consider $(d_{1}, d_{2}) \geq (1, 1)$. If $\operatorname{Hyp}_{(d_{1}, d_{2})}(B, 2) \neq \emptyset$, then: 
    \begin{enumerate}
        \item\label{prop: F0 Hyp(B,2)=E(B).(1)} $B_{1} \in \lvert (1, 1) \rvert$.
        \item\label{prop: F0 Hyp(B,2)=E(B).(2)} $d_{1} = d_{2} = 1$.
        \item\label{prop: F0 Hyp(B,2)=E(B).(3)}
        Furthermore, let $\mathfrak{I} \subseteq \{1,2,3\}$ be the set of indices such that $B_i \in \lvert (1,1) \rvert$ for all $i \in \mathfrak{I}$. Then
        $$
        \operatorname{Hyp}_{(1, 1)}(B, 2) = \bigcup\limits_{\substack{ 
        i\in\mathfrak{I} \\ 
        \{i,j,k\} = \{1,2,3\} \\
        p_{i,j}\in B_i \cap B_j \\
        p_{i,k}\in B_i \cap B_k}} 
        \operatorname{Hyp}_{(1, 1)}(B_j, p_{i, j}) \cap \operatorname{Hyp}_{(1, 1)}(B_k, p_{i, k}).
        $$
    \end{enumerate}
    In particular,
    $$
    \operatorname{Hyp}(B, 2) 
    = \operatorname{Hyp}_{(1, 0)}(B, 2) \cup \operatorname{Hyp}_{(0, 1)}(B, 2) \cup \operatorname{Hyp}_{(1, 1)}(B, 2) 
    = \mathcal{E}(B).
    $$
\end{proposition}

\begin{proof}
    Let $D \in \operatorname{Hyp}_{(d_{1}, d_{2})}(B, 2)$. Then, $\#(D \cap B) \leq 2$, and $(D \cdot B_{i}) \geq 2 > 0$ for any $i = 1, 2, 3$. On the other hand, $B_{1} \cap B_{2} \cap B_{3} = \emptyset$ by assumption. Therefore, we must have $\#(D \cap B) = 2$. This implies that any point in $D\cap B$ is a unibranch point on $D$. 
    
    \begin{claim}
        $D \cap B_{3} \subseteq N = \bigcup_{i<j} B_i \cap B_j$, and $\#(D \cap B_{3}) = 1$.
    \end{claim}
    \begin{proof}[Proof of the Claim]
        Suppose there exists a point $q$ in $D \cap B_{3}$ that is not contained in $N$. Then, $D$ must intersect $B_{1} \cup B_{2}$ in exactly one point. Since $(B_{i} \cdot B_{j}) = \alpha_i \beta_j + \alpha_j \beta_i \geq 2$, we know $B_{i} \cap B_{j} \neq \emptyset$ for any $i \neq j$. Thus, $D$ intersects $B_{1} \cup B_{2}$ at a unique point $p_{1,2} \in B_{1} \cap B_{2}$, and it must be a unibranch $n$-fold point of $D$ for some $n \geq 1$. By Lem.~\ref{lem: Upper bound of unibranch multiplicity on Hirzebruch Surfaces}, we obtain
        $
        n \leq \min\{d_{1}, d_{2}\} < d_{1} + d_{2}
        $.
        
        As $B_{1}$ and $B_{2}$ meet transversally at $p_{1, 2}$, we see that $D$ is tangent to one of them and transverse to the other. Assume $D$ is transverse to $B_{i}$, $i \in \{1, 2\}$, then
        $$
        d_{1} \beta_{i} + d_{2} \alpha_{i}
        = (D \cdot B_{i})
        = (D \cdot B_{i})_{p_{1, 2}}
        = n
        \leq \min\{d_{1}, d_{2}\}
        < d_{1} + d_{2},
        $$
        which is a contradiction. Thus, we have proved that $D \cap B_{3} \subseteq N$.

        Now, if $\#(D \cap B_{3}) = 2$, then $D \cap B_{3} = \{p_{1, 3}, p_{2, 3}\}$ for some points \( p_{1,3} \in B_{1}\cap B_{3} \) and \( p_{2,3} \in B_{2}\cap B_{3} \). As $B_{3}$ meets $B_{i}$ transversally at $p_{i, 3}$ for any $i = 1, 2$, we see that $D$ must be transverse to $B_{i}$ or $B_{3}$ at $p_{i, 3}$.

        If $D$ is transverse to $B_{i}$ at $p_{i, 3}$ for some $i \in \{1,2\}$, as before, we may assume $p_{i, 3}$ is a unibranch $n$-fold point of $D$ for some $n \geq 1$. Then
        $$
        d_{1} \beta_{i} + d_{2} \alpha_{i}
        = (D \cdot B_{i})
        = (D \cdot B_{i})_{p_{i, 3}}
        = n
        \leq \min\{d_{1}, d_{2}\}
        < d_{1} + d_{2},
        $$
        which is again a contradiction. Hence, $D$ is transverse to $B_{3}$ at both $p_{1, 3}$ and $p_{2, 3}$. Let $m_{i} \coloneqq \operatorname{mult}_{p_{i, 3}}(D)$ for $i=1,2$, then
        $
        m_{1} + m_{2}
        \leq 2 \min\{d_{1}, d_{2}\}
        \leq d_{1} + d_{2}
        $.
        Therefore, we have
        $$
        d_{1} \beta_{3} + d_{2} \alpha_{3}
        = (D \cdot B_{3})
        = (D \cdot B_{3})_{p_{1, 3}} + (D \cdot B_{3})_{p_{2, 3}}
        = m_{1} + m_{2}
        \leq d_{1} + d_{2},
        $$
        which is possible only if $\alpha_{3} = \beta_{3} = 1$. In this case, we must have
        $
        \#(D \cap B_{1}) = \#(D \cap B_{2}) = 1
        $.
        And by the assumption $\alpha_1 + \beta_1 \leq \alpha_2 + \beta_2 \leq \alpha_3 + \beta_3$, we also have
        $
        B_{1}, B_{2}, B_{3} \in \lvert (1,1) \rvert
        $.
        As $D \cap B_{1} = \{p_{1, 3}\} \subset N$, we see that $B_{1}$ satisfies the requirement of the claim if it is in the place of $B_{3}$. Since in this case, we have $\alpha_i = \beta_i =1$ for all $i\in \{1,2,3\}$, we may re-index the three irreducible components such that $B_1$ and $B_3$ switch their places. This completes the proof of the claim. 
    \end{proof}

    Now, we have proved that $D \cap B_{3} = \{ p_{i, 3} \}$ for some $i \in \{1, 2\}$. Notice that $D$ must be tangent to $B_{3}$ at $p_{i,3}$ because
    $
    (D \cdot B_{3}) = \alpha_{3} d_{1} + \beta_{3} d_{2} \geq d_{1} + d_{2} > \operatorname{mult}_{p_{i,3}}(D)
    $.
    
    Thus, $D$ must be transverse to $B_{i}$ at $p_{i, 3}$. Then $D$ must meet $B_{i}$ at a further point, since
    $$
    (D \cdot B_{i}) 
    \geq d_1 + d_2
    > \operatorname{mult}_{p_{i,3}}
    = (D \cdot B_i)_{p_{i,3}}.
    $$
    As $D$ must meet the other component $B_{j}$ (where $j \neq i, 3$), we know that there exists a point $p_{1, 2} \in B_{1} \cap B_{2}$ such that $D \cap B_{j} = \{p_{1,2}\}$ and $D \cap B = \{p_{1,2}, p_{i,3}\}$. 

    Now, $B_{j}$ and $D$ meet only at $p_{1, 2}$. Since $(B_{j} \cdot D) \geq d_1 + d_2 > \operatorname{mult}_{p_{1,2}}(D)$, we see that $D$ cannot be transverse to $B_{j}$ at $p_{1, 2}$. Hence, it must be transverse to $B_{i}$ at $p_{1, 2}$. By Lem.~\ref{lem: Upper bound of unibranch multiplicity on Hirzebruch Surfaces}, we know 
    $
    \operatorname{mult}_{p_{1, 2}}(D) \leq \min\{d_{1}, d_{2}\}
    $
    and
    $
    \operatorname{mult}_{p_{i, 3}}(D) \leq \min\{d_{1}, d_{2}\}.
    $
    Then
    \begin{multline}
        d_{1} \beta_{i} + d_{2} \alpha_{i}
        = (D \cdot B_{i})
        = (D \cdot B_{i})_{p_{1, 2}} + (D \cdot B_{i})_{p_{i, 3}} \\
        = \operatorname{mult}_{p_{1, 2}}(D) + \operatorname{mult}_{p_{i, 3}}(D) 
        \leq 2 \min\{d_{1}, d_{2}\}
        \leq d_{1} + d_{2}.
    \end{multline}
    Hence the equalities must hold, we obtain $\alpha_{i} = \beta_{i} = 1$ and $d_1 = d_2 = \operatorname{mult}_{p_{1, 2}}(D) = \operatorname{mult}_{p_{i, 3}}(D)$. By the assumption $\alpha_1 + \beta_1 \leq \alpha_2 + \beta_2 \leq \alpha_3 + \beta_3$, we get $B_{1} \in \lvert (1,1) \rvert$, which proves part (\ref{prop: F0 Hyp(B,2)=E(B).(1)}) of the proposition.

    Without loss of generality, let $i=1$. Then $D \cap B_3 = \{ p_{1,3} \}$. Denote $d \coloneqq d_{1} = d_{2}$. Then, $D \in \lvert (d, d) \rvert$ for some $d \geq 1$, and 
    $$
    D \in \operatorname{Hyp}_{(d, d)}(B, 2)
    \subseteq \operatorname{Hyp}_{(d, d)}^{d}(B_{2}, p_{1, 2}) \cap \operatorname{Hyp}_{(d, d)}^{d}(B_{3}, p_{1, 3}).
    $$
    Thus, both $p_{1, 2}$ and $p_{1, 3}$ are unibranch $d$-fold points of $D$. As discussed above, $D$ meets $B_{1}$ transversally at $p_{1, 2}$ and $p_{1, 3}$. Therefore, $D$ must be hypertangent to $B_{j}$ at $p_{1, j}$ for all $j = 2, 3$. Then, by Thm.~\ref{thm: bounddeltainvariant}, for all $j=2,3$, we obtain 
    $$
    \delta_{D}(p_{1, j}) 
    \geq \frac{(d - 1)(d (\alpha_{j} + \beta_{j}) - 1) + d - 1}{2} 
    \geq \frac{(d - 1)(2 d)}{2}
    = d (d - 1).
    $$
    Then
    \begin{multline}
        0 \leq p_{a}(D^{\nu})
            \leq p_{a}(D) - \delta_{D}(p_{1, 3}) - \delta_{D}(p_{2, 3}) \\
            \leq (d - 1)^{2} - d(d - 1) - d(d - 1) 
            = -(d + 1)(d - 1)
            \leq 0.
    \end{multline}
    Therefore, the equalities must hold, which implies that $d = 1$ and $p_{a}(D^{\nu}) = 0$. This proves part (\ref{prop: F0 Hyp(B,2)=E(B).(2)}) of the proposition.
 
    Following the notation of part (\ref{prop: F0 Hyp(B,2)=E(B).(3)}) of the proposition, if 
    $\operatorname{Hyp}_{(1,1)}(B,2)
    \neq \emptyset$, then we must have $1 \in \mathfrak{I}$. If $\mathfrak{I}=\{1\}$, i.e., 
    $\{ B_1 , B_2 , B_3\} \cap \lvert (1,1) \rvert = \{ B_1 \}$, then we must have $i=1$. In this case, we have proved 
    $$
    \operatorname{Hyp}_{(1, 1)}(B, 2)
    \subseteq \bigcup\limits_{\substack{p_{1,2}\in B_1 \cap B_2 \\ p_{1,3}\in B_1 \cap B_3}} \operatorname{Hyp}_{(1, 1)}(B_{2}, p_{1, 2}) \cap \operatorname{Hyp}_{(1, 1)}(B_{3}, p_{1, 3}),
    $$
    while the reverse inclusion 
    $$
    \operatorname{Hyp}_{(1, 1)}(B, 2)
    \supseteq  \bigcup\limits_{\substack{p_{1,2}\in B_1 \cap B_2 \\ p_{1,3}\in B_1 \cap B_3}} \operatorname{Hyp}_{(1, 1)}(B_{2}, p_{1, 2}) \cap \operatorname{Hyp}_{(1, 1)}(B_{3}, p_{1, 3})
    $$
    is obvious. 

    If $B_2$ or $B_3$ is also contained in $\lvert (1,1) \rvert$, then we may switch it with $B_1$, and the preceding arguments remain valid. By symmetry, we immediately obtain part (\ref{prop: F0 Hyp(B,2)=E(B).(3)}) of the proposition.

    Hence, we conclude that
    $$
    \operatorname{Hyp}(B, 2) 
    = \operatorname{Hyp}_{(1, 0)}(B, 2) \cup \operatorname{Hyp}_{(0, 1)}(B, 2) \cup \operatorname{Hyp}_{(1, 1)}(B, 2) 
    = \mathcal{E}(B).
    $$ 
    This completes the proof of the theorem.
\end{proof}

\smallskip

Now, we give an effective bound for $\mathcal{E}(B)$ in Prop.~\ref{prop: F0 Hyp(B,2)=E(B)}:

\begin{theorem}\label{thm: F0 bound E(B)}
    Consider $\mathbb{F}_0$. Let $B$ be as in Prop.~\ref{prop: F0 Hyp(B,2)=E(B)}. Then $\mathcal{E}(B)$ is finite. More precisely, 
    $$
    \lvert \mathcal{E}(B) \rvert 
    \leq 
    \begin{cases}
        2 \lvert N \rvert 
        & \text{ if } \{ B_1 , B_2 , B_3 \} \cap \lvert (1,1) \rvert = \emptyset;\\ 
        (\alpha_{2} + \beta_{2})(\alpha_{3} + \beta_{3}) + 2 \lvert N \rvert
        & \text{ if } \{ B_1 , B_2 , B_3 \} \cap \lvert (1,1) \rvert = \{ B_1 \}; \\
        4(\alpha_{3} + \beta_{3}) + 2 \lvert N \rvert
        & \text{ if } \{ B_1 , B_2 , B_3 \} \cap \lvert (1,1) \rvert = \{ B_1 , B_2 \}; \\
        12 + 2 \lvert N \rvert = 24 
        & \text{ if } \{ B_1 , B_2 , B_3 \} \subset \lvert (1,1) \rvert.
    \end{cases}
    $$
    Here 
    $\lvert N \rvert
    = \sum_{i < j} (B_{i} \cdot B_{j})
    = \sum_{i < j} (\alpha_{i} \beta_{j} + \alpha_{j} \beta_{i})$.
\end{theorem}

\begin{proof}
    If $D$ is an integral curve in $\operatorname{Hyp}_{(1, 0)}(B, 2)$, then $(D \cap B_{i}) \geq 1 >0$ for any $i = 1, 2, 3$, since $(\alpha_i , \beta_i) \geq (1, 1)$ by assumption. Hence, $D$ must contain some point in $N = \bigcup_{i < j} B_{i} \cap B_{j}$. Since for any point on $\mathbb{F}_{0} \cong \mathbb{P}^{1} \times \mathbb{P}^{1}$, there is precisely one integral curve in the linear system $\lvert (1, 0) \rvert$ passing through it, it follows that
    $
    \lvert \operatorname{Hyp}_{(1,0)}(B,2) \rvert
    \leq \lvert N \rvert
    $.
    By symmetry, we obtain
    $
    \lvert \operatorname{Hyp}_{(0,1)}(B, 2) \rvert
    \leq \lvert N \rvert
    $.

    Recall that by Prop.~\ref{prop: F0 Hyp(B,2)=E(B)}, we have
    $$
    \mathcal{E}(B) 
    = \operatorname{Hyp} (B,2)
    = \operatorname{Hyp}_{(1,0)} (B,2) \cup \operatorname{Hyp}_{(0,1)} (B,2) \cup \operatorname{Hyp}_{(1,1)} (B,2).
    $$

    If $\{ B_1 , B_2 , B_3\} \cap \lvert (1,1) \rvert = \emptyset$, then, by Prop.~\ref{prop: F0 Hyp(B,2)=E(B)}, we know 
    $\operatorname{Hyp}_{(1,1)}(B, 2) = \emptyset$ 
    and 
    $\mathcal{E}(B) 
    = \operatorname{Hyp}_{(1,0)}(B, 2) \cup \operatorname{Hyp}_{(0,1)}(B, 2)$. In this case, we obtain 
    $\lvert \mathcal{E}(B) \rvert \leq 2\lvert N \rvert$.
   
    If $\operatorname{Hyp}_{(1,1)}(B, 2) \neq \emptyset$, then we must have $B_{1} \in \lvert (1,1) \rvert$. 
    
    Let $D$ be a curve in $\operatorname{Hyp}_{(1,1)}(B,2)$.
    
    If $\{ B_1 , B_2, B_3 \} \cap \lvert (1,1) \rvert = \{ B_1 \}$, then, as in the proof of Prop.~\ref{prop: F0 Hyp(B,2)=E(B)}, we know 
    $D \cap B_1 = \{p_{1,2} , p_{1,3}\}$ for some points 
    $p_{1,2} \in B_1 \cap B_2$
    and
    $p_{1,3} \in B_1 \cap B_3$. 
    With $p_{1,2}$ and $p_{1,3}$ fixed, there is at most one curve in $\operatorname{Hyp}_{(1,1)}(B, 2)$ that contains both $p_{1,2}$ and $p_{1,3}$. Indeed, assume there are two such curves $D_{1}, D_{2} \in \operatorname{Hyp}_{(1,1)}(B, 2)$, then we have
    $$
    D_{1} \cap B_{2} = \{p_{1,2}\} = D_{2} \cap B_{2}
    \quad
    \text{and}
    \quad
    D_{1} \cap B_{3} = \{p_{1,3}\} = D_{2} \cap B_{3}.
    $$
    Then
    $
    (D_{1} \cdot B_{2})_{p_{1, 2}} = (D_{1} \cdot B_{2}) = \alpha_{2} + \beta_{2} = (D_{2} \cdot B_{2}) = (D_{2} \cdot B_{2})_{p_{1, 2}}
    $.
    
    Since $(D_{1} \cdot D_{2}) = 2$ and $D_{1} \cap D_{2} = \{p_{1,2}, p_{1, 3}\}$, we must have
    $
    (D_{1} \cdot D_{2})_{p_{1,2}} = 1 = (D_{1} \cdot D_{2})_{p_{1,3}}
    $.
    This leads to
    $
    (D_{1} \cdot D_{2})_{p_{1, 2}} 
    = 1
    < \alpha_{2}+\beta_{2}
    = (D_{1} \cdot B_{2})_{p_{1,2}} 
    = (D_{2} \cdot B_{2})_{p_{1,2}}
    $.
    Note that $p_{1, 2}$ is a smooth point on all three curves $D_{1}, D_{2}$, and $B_{2}$. But then, by Thm.~\ref{thm: Strong Triangle Inequality}, the smallest two among
    $$
    (D_{1} \cdot D_{2})_{p_{1,2}},
    \quad
    (D_{1} \cdot B_{2})_{p_{1,2}},
    \quad
    (D_{2} \cdot B_{2})_{p_{1,2}}
    $$
    are equal---a contradiction.
    
    Since there are $(B_{1} \cdot B_{2})$ choices for $p_{1, 2}$ and $(B_{1} \cdot B_{3})$ choices for $p_{1, 3}$, we conclude that 
    $$
    \lvert \operatorname{Hyp}_{(1,1)}(B, 2) \rvert 
    \leq (B_{1} \cdot B_{2}) (B_{1} \cdot B_{3})
    = (\alpha_2 + \beta_2) (\alpha_3 + \beta_3).
    $$
    Summing up, we obtain
    $$
    \lvert \mathcal{E}(B) \rvert 
    \leq (\alpha_{2} + \beta_{2})(\alpha_{3} + \beta_{3}) + 2 \lvert N \rvert
    $$
    in the case 
    $\{ B_1 , B_2, B_3 \} \cap \lvert (1,1) \rvert = \{ B_1 \}$ 
    as desired.

    If $\{ B_1 , B_2, B_3 \} \cap \lvert (1,1) \rvert = \{ B_1 , B_2 \}$, then we may switch $B_1$ and $B_2$. By symmetry, we obtain
    $$
    \lvert \operatorname{Hyp}_{(1,1)}(B,2) \rvert
    \leq (B_1 \cdot B_2)(B_1 \cdot B_3) + (B_2 \cdot B_1)(B_2 \cdot B_3)
    = 4(\alpha_3 + \beta_3).
    $$
    Hence,
    $$
    \lvert \mathcal{E}(B) \rvert
    \leq 4(\alpha_3 + \beta_3) + 2 \lvert N \rvert.
    $$

    If $\{ B_1 , B_2, B_3 \} \subset \lvert (1,1) \rvert$, then we may switch $B_2$ and $B_3$ with $B_1$. By symmetry, we obtain
    $$
    \lvert \operatorname{Hyp}_{(1,1)}(B,2) \rvert
    \leq (B_1 \cdot B_2)(B_1 \cdot B_3) + (B_2 \cdot B_1)(B_2 \cdot B_3) + (B_3 \cdot B_1)(B_3 \cdot B_2)
    = 12.
    $$
    Note that $\lvert N \rvert = 6$ in this case, hence
    $$
    \lvert \mathcal{E}(B) \rvert
    \leq 12 + 2 \lvert N \rvert
    = 24.
    $$
    
    This completes the proof.
\end{proof}

\begin{proposition}\label{prop: emptyness of exceptional set F0}
Consider $\mathbb{F}_0$. Let $B$ be as in Prop.~\ref{prop: F0 Hyp(B,2)=E(B)}. 

If 
$\alpha_i \geq 3$ and $\beta_i \geq 3$ for all $i\in \{1,2,3\}$ and $B$ is general, then
$
\mathcal{E}(B) = \emptyset
$.
\end{proposition}
\begin{proof}
    By Prop.~\ref{prop: F0 Hyp(B,2)=E(B)}, we know $\operatorname{Hyp}(B,2) = \operatorname{Hyp}_{(1,0)}(B,2) \cup \operatorname{Hyp}_{(0,1)}(B,2)$.
    
    We first note that $(1,0)\vert_{B_i}$ is a base point free linear system on $B_i$ for any $i \in \{1,2,3\}$. For any $p \in B_i$, let 
    $F_p^{(1,0)}$ be the unique curve in $ \lvert (1,0) \rvert$ that contains $p$. If $B_i$ is smooth, then by Bertini's theorem (\cite[Ch.~\Romannum{3}, Cor.~10.9, p.~274]{Hartshorne1977}), we know for a general point $p \in B_i$ that $F_p^{(1,0)}$ intersects $B_i$ transversally. By symmetry, we can define $F_p^{(0,1)}$, and the same holds. 
    Define the following set of points:
    $$
    T_f(B_i)
    \coloneqq
    \{
    p \in B_i
    \
    \vert
    \ 
    F_p^{(1,0)} \text{ or } F_p^{(0,1)}
    \text{ is tangent to } B_i
    \}.
    $$
    Then $T_f(B_i)$ is a finite set.

    We choose $B$ general in the following sense:
    $B_1 \in \lvert (\alpha_1 , \beta_1) \lvert$ is smooth; 
    $B_2 \in \lvert (\alpha_2 , \beta_2) \rvert$ is smooth, 
    intersects $B_1$ transversally, 
    and 
    $B_2 \cap T_f(B_1) = \emptyset$;
    and $B_3 \in \lvert (\alpha_3 , \beta_3) \rvert$ is smooth,
    intersects $B_1 \cup B_2$ transversally, and 
    $B_3 \cap T_f(B_1) = \emptyset = B_3 \cap T_f(B_2)$.

    Then, by construction, 
    for any $i < j$, $i,j \in \{1,2,3\}$ and for any point $p_{i,j} \in B_i \cap B_j$, $F_{p_{i,j}}^{(1,0)}$ intersects $B_i$ transversally. We obtain
    $$
    \# F_{p_{i,j}}^{(1,0)} \cap B
    \geq 
    \# F_{p_{i,j}}^{(1,0)} \cap B_i
    = (F_{p_{i,j}}^{(1,0)} \cdot B_i)
    = \beta_i
    \geq 3.
    $$
    Thus $F_{p_{i,j}}^{(1,0)} \not\in \operatorname{Hyp}_{(1,0)} (B,2)$. By the first part of the proof of Thm.~$\ref{thm: F0 bound E(B)}$, we know 
    $
    \operatorname{Hyp}_{(1,0)} (B,2) = \emptyset.
    $
    By symmetry, we also obtain
    $
    \operatorname{Hyp}_{(0,1)} (B,2) = \emptyset.
    $
    Combining with Prop.~\ref{prop: F0 Hyp(B,2)=E(B)}, we conclude that
    $$
    \mathcal{E}(B) 
    = \operatorname{Hyp} (B,2)
    = \operatorname{Hyp}_{(1,0)} (B,2) \cup \operatorname{Hyp}_{(0,1)} (B,2)
    = \emptyset.
    $$
\end{proof}

Now, we consider a $3C$-curve $B = B_1 \cup B_2 \cup B_3$ on the Hirzebruch surface $\mathbb{F}_{e}$ with $e \geq 2$ such that $K_{\mathbb{F}_{e}} + \mathcal{L}(B)$ is big. 
Assume $B \in \lvert \alpha C_{1} + \beta f \rvert$ and $B_i \in \lvert \alpha_i C_1 + \beta_i f\rvert$ for all $i = 1,2,3$. 

We note that $\beta \geq -e$. Indeed, by Cor.~\ref{cor: Irreducible Curves on Hirzebruch Surfaces}, if, say, $\beta_1 < 0$, then $\beta_1 = -e$, $\alpha_1 = 1$ and $B_1 = C_0$. Since $B$ is reduced and contains at most one copy of $C_{0}$ as its irreducible component, we get $\beta_2 , \beta_3 \geq 0$, and it follows that $\beta \geq -e$.

Since $K_{\mathbb{F}_{e}} = -2 C_{1} + (e - 2) f$, by Prop.~\ref{prop: Criterion for bigness on Hirzebruch Surfaces}, we have that $K_{\mathbb{F}_{e}} + \mathcal{L}(B)$ is big if and only if $\alpha > 2$ and $(\alpha - 2) e + (\beta + e - 2) > 0$. The latter condition is equivalent to $\beta + (\alpha - 1) e > 2$.

If $B = B_{1} \cup B_{2} \cup B_{3}$ is a $3C$-curve on $\mathbb{F}_{e}, e \geq 2$, where each $B_{i} \in \lvert \alpha_{i} C_{1} + \beta_{i} f \rvert$ satisfies $\alpha_{i} \geq 1$ and $\beta_{i} \geq 0$ for $i = 1, 2, 3$, then $\alpha \geq 3$ and $\beta \geq 0$. Consequently, $K_{\mathbb{F}_{e}} + \mathcal{L}(B)$ is big. We have the following:

\begin{proposition}\label{prop: Fe Hyp(B,2)=E(B)}
    Consider $\mathbb{F}_{e}, e\geq 2$. Let $B=B_{1}\cup B_{2}\cup B_{3}  \in \lvert \alpha C_{1} + \beta f \rvert$ be a $3C$-curve on $\mathbb{F}_{e}$, and denote $B_{i}\in \lvert \alpha_{i}C_{1}+\beta_{i}f \rvert$ for all $i=1,2,3$. Assume that $K_{\mathbb{F}_e} + \mathcal{L}(B)$ is big, and that none of $B_{1}, B_{2}, B_{3}$ belongs to $\lvert f \rvert$ or is $C_{0}$. 

    Then
    $$
    \operatorname{Hyp}(B,2)
    = \operatorname{Hyp}_{C_{1}}(B,2) \cup \operatorname{Hyp}_{C_{0}}(B,2) \cup \operatorname{Hyp}_{f}(B,2) 
    = \mathcal{E}(B),
    $$
    and $\operatorname{Hyp}_{C_{1}}(B,2) \neq \emptyset$ only if $e = 2$ and one of the three components lies in $\lvert C_{1} \rvert$.
\end{proposition}

\begin{proof}
    By assumption, we have $\alpha_{i}\geq 1$ and $\beta_{i}\geq 0$ for all $i=1,2,3$. Let $D \in \lvert d_{1} C_{1} + d_{2} f \rvert$ be an integral curve such that $d_{1} > 0$ and $d_{2} \geq 0$, which is equivalent to $D$ being neither contained in $\lvert f \rvert$ nor equal to $C_{0}$. Then $(D \cdot B_i) \geq e \geq 2 > 0$ for all $i \in \{1,2,3\}$, hence $D\cap B_i \neq \emptyset$.
    
    Recall that $N = \bigcup_{i < j} B_{i} \cap B_{j}$. If $D \in \operatorname{Hyp}(B,2)$, then $\#(D \cap B) \leq 2$. Since $B_{1} \cap B_{2} \cap B_{3} = \emptyset$, we have $\#(D \cap B) = 2$. This implies that any point in $D\cap B$ is a unibranch point on $D$.

    \begin{claim}
        $D \cap B_{i} \subseteq N$ for all $i \in \{1, 2, 3\}$.
    \end{claim}
    \begin{proof}[Proof of the Claim]
        Without loss of generality, we may assume $i=1$. Suppose there exists a point $q$ in $D \cap B_{1}$ that is not contained in $N$, then $D$ must intersect $B_{2} \cup B_{3}$ in exactly one point $p_{2,3} \in B_2 \cap B_3$.

        Therefore, $D$ is tangent to one of $B_{2}$ and $B_{3}$ at $p_{2,3}$ and transverse to the other. Without loss of generality, assume $D$ is transverse to $B_{2}$ and $p_{2,3}$ is a unibranch $n$-fold point of $D$. Then
        $$
        n = (D \cdot B_{2}) = \left( (d_{1} C_{1} + d_{2} f) \cdot (\alpha_{2} C_{1} + \beta_{2} f) \right) = e \alpha_{2} d_{1} + \alpha_{2} d_{2} + \beta_{2} d_{1} \geq e \alpha_{2} d_{1} > d_{1},
        $$
        which contradicts Lem.~\ref{lem: Upper bound of unibranch multiplicity on Hirzebruch Surfaces}. Thus, we have $D \cap B_{1} \subseteq N$ as desired.
    \end{proof}    

    Now, since $D \cap B_{i} \subset N$ for any $i = 1, 2, 3$, we know there exists $j$ such that $\#(D \cap B_{j}) = 2$ and $\#(D \cap B_{i}) = 1$ for any $i \neq j$. Without loss of generality, let $j = 1$, that is, $\#(D \cap B_{1}) = 2$. Then 
    $
    D \cap B = D \cap B_{1} = \{p_{1,2}, p_{1,3}\}
    $,
    where $p_{1,2} \in B_{1} \cap B_{2}$ and $p_{1,3} \in B_{1} \cap B_{3}$.

    If $D$ is transverse to $B_{i}$ at $p_{1,i}$ for some $i \in \{2,3\}$, as before, assume $p_{1,i}$ is a unibranch $n$-fold point of $D$ for some $n \geq 1$. Then
    \begin{multline}
        n 
        = (D \cdot B_{i})_{p_{1,i}} 
        = (D \cdot B_{i}) = \left( (d_{1} C_{1} + d_{2} f) \cdot (\alpha_{i} C_{1} + \beta_{i} f) \right) 
        \\
        = e \alpha_{i} d_{1} + \alpha_{i} d_{2} + \beta_{i} d_{1} 
        \geq e d_{1} 
        > d_{1},
    \end{multline}
    which contradicts Lem.~\ref{lem: Upper bound of unibranch multiplicity on Hirzebruch Surfaces}.

    Hence, $D$ is transverse to $B_{1}$ at both $p_{1,2}$ and $p_{1,3}$. Let $m_{i} \coloneqq \operatorname{mult}_{p_{1,i}}(D)$ for $i = 2, 3$. Then
    \begin{multline}
    2 d_{1} 
    \geq m_{2} + m_{3} 
    = (D \cdot B_{1})_{p_{1,2}} + (D \cdot B_{1})_{p_{1,3}} 
    = (D \cdot B_{1}) 
    \\
    = e \alpha_{1} d_{1} + \alpha_{1} d_{2} + \beta_{1} d_{1} 
    \geq e d_{1} 
    \geq 2 d_{1}.    
    \end{multline}
    Thus, all the equalities must hold, which implies
    $$
    e = 2, \quad \alpha_{1} = 1, \quad \beta_{1} = 0, \quad d_{2} = 0, \quad m_{2} = m_{3} = d_{1}.
    $$

    From this, we conclude that 
    $$
    \operatorname{Hyp}_{d_{1} C_{1} + d_{2} f}(B,2) = \emptyset
    $$
    if
    $e \geq 3$, $d_{1} > 0$, $d_{2} \geq 0$
    or
    $e = 2$, $d_{1} > 0$, $d_{2} > 0$.

    We are left with the case $e = 2$, $\alpha_{1} = 1$, $\beta_{1} = 0$, $d_{2} = 0$, and $m_{2} = m_{3} = d_{1}$. Now, $D \in \lvert d_{1} C_{1} \rvert$. By Thm.~\ref{thm: bounddeltainvariant}, for any $i=2,3$, we have
    $$
    \delta_{D}(p_{1,i}) \geq \frac{(d_{1} - 1) \left( (D \cdot B_{i}) - 1 \right)}{2} = \frac{(d_{1} - 1)(2 \alpha_{i} d_{1} + \beta_{i} d_{1} - 1)}{2}.
    $$
    Hence, we obtain
    \begin{equation}
        \begin{split}
            0 \leq p_{a}(D^{\nu}) & \leq p_{a}(D) - \delta_{D}(p_{1,2}) - \delta_{D}(p_{1,3}) \\
            & \leq \frac{1}{2} (d_{1}-1) (2d_{1}-2) 
            - \frac{1}{2} (d_{1}-1) (2\alpha_{2}d_{1}+\beta_{2}d_{1}-1)\\
            & \quad - \frac{1}{2} (d_{1}-1) (2\alpha_{3}d_{1}+\beta_{3}d_{1}-1) \\
            & = \frac{1}{2}(d_{1}-1)d_{1}(2-2\alpha_{2}-2\alpha_{3}-\beta_{2}-\beta_{3})
            = -(d_{1}-1)d_{1}
            \leq 0
            .
        \end{split}
    \end{equation}
    Thus, all the equalities here must hold, which implies $d_{1} = 1$ and $p_{a}(D^{\nu}) = 0$. Therefore, when $e = 2$, for $d_{1} > 0$ and $d_{2} \geq 0$, the set $\operatorname{Hyp}_{d_{1} C_{1} + d_{2} f}(B,2)$ is non-empty only if $d_{1} = 1$, $d_{2} = 0$, and $B_1 \in \lvert C_1 \rvert$. This completes the proof.
\end{proof}

We obtain an effective bound for $\mathcal{E}(B)$ in Prop.~\ref{prop: Fe Hyp(B,2)=E(B)}:

\begin{theorem}\label{thm: Fe bound E(B)}
    Consider $\mathbb{F}_e , e\geq 2$. Let $B$ be as in Prop.~\ref{prop: Fe Hyp(B,2)=E(B)}. Then $\mathcal{E}(B)$ is finite. More precisely:
    
    If $e \geq 3$, then
    $$
    \lvert \mathcal{E}(B) \rvert 
    \leq 
    1 + \lvert N \rvert.
    $$
    If $e=2$, assume 
    $\alpha_1 + \beta_1 \leq \alpha_2 + \beta_2 \leq \alpha_3 + \beta_3$, then
    $$
    \lvert \mathcal{E}(B) \rvert 
    \leq
    \begin{cases}
        1 + \lvert N \rvert 
        & \text{ if } \{ B_1 , B_2 , B_3 \} \cap \lvert C_1 \rvert = \emptyset;\\ 
        1 + (2 \alpha_{2} + \beta_{2})(2 \alpha_{3} + \beta_{3}) + \lvert N \rvert
        & \text{ if } \{ B_1 , B_2 , B_3 \} \cap \lvert C_1 \rvert = \{ B_1 \}; \\
        1+4(2\alpha_{3} + \beta_{3}) + \lvert N \rvert
        & \text{ if } \{ B_1 , B_2 , B_3 \} \cap \lvert C_1 \rvert = \{ B_1 , B_2 \}; \\
        1+ 12 + \lvert N \rvert = 19 
        & \text{ if } \{ B_1 , B_2 , B_3 \} \subset \lvert C_1 \rvert.
    \end{cases}
    $$
    Here 
    $\lvert N \rvert
    = \sum_{i < j} (B_{i} \cdot B_{j})
    = \sum_{i < j} (e\alpha_i \alpha_j + \alpha_{i} \beta_{j} + \alpha_{j} \beta_{i})$.
\end{theorem}

\begin{proof}
    First, by Rem.~\ref{rmk: C0 unique negative self intersection curve on Fe}, we have $\lvert \operatorname{Hyp}_{C_{0}}(B,2) \rvert \leq 1$.

    Next, by the same argument as in the proof of Thm.~\ref{thm: F0 bound E(B)} for the $\mathbb{F}_{0}$ case, we see that any integral curve $D \in \operatorname{Hyp}_{f}(B,2)$ contains some point in $N$. Therefore, we obtain
    $$
    \lvert \operatorname{Hyp}_{f}(B,2) \rvert
    \leq \sum_{i < j} (B_{i} \cdot B_{j})
    = \lvert N \rvert.
    $$

    Recall that by Prop.~\ref{prop: Fe Hyp(B,2)=E(B)}, we have
    $$
    \mathcal{E}(B)
    = \operatorname{Hyp} (B,2)
    = \operatorname{Hyp}_{C_0} (B,2) \cup 
    \operatorname{Hyp}_{C_1} (B,2) \cup \operatorname{Hyp}_{f} (B,2).
    $$

    If $e \geq 3$, then $\operatorname{Hyp}_{C_1}(B,2) = \emptyset$ and 
    $\mathcal{E}(B) = \operatorname{Hyp}_{C_0}(B,2) \cup \operatorname{Hyp}_{f}(B,2)$. Hence
    $$
    \lvert \mathcal{E}(B) \rvert 
    \leq 
    1 + \sum_{i < j} (B_{i} \cdot B_{j}) 
    = 1 + \lvert N \rvert.
    $$

    If $e = 2$ and $\{B_1 , B_2 , B_3\} \cap \lvert C_1 \rvert = \emptyset$, then by Prop.~\ref{prop: Fe Hyp(B,2)=E(B)}, we know $\operatorname{Hyp}_{C_{1}}(B,2) = \emptyset$ and 
    $\mathcal{E}(B) = \operatorname{Hyp}_{C_0}(B,2) \cup \operatorname{Hyp}_{f}(B,2)$. Hence 
    $$
    \lvert \mathcal{E}(B) \rvert \leq 1 + \lvert N \rvert.
    $$

    If $e = 2$ and $\operatorname{Hyp}_{C_{1}}(B,2) \neq \emptyset$, then $\{B_1 , B_2 , B_3\} \cap \lvert C_1 \rvert \neq \emptyset$. By the assumption 
    $\alpha_1 + \beta_1 
    \leq \alpha_2 + \beta_2  
    \leq \alpha_3 + \beta_3$, 
    we must have $B_1 \in \lvert C_1 \rvert$.

    Let $D$ be a curve in $\operatorname{Hyp}_{C_1}(B,2)$.

    If $\{B_1 , B_2 , B_3\} \cap \lvert C_1 \rvert = \{B_1\}$, then, as in the proof of Prop.~\ref{prop: Fe Hyp(B,2)=E(B)}, we know 
    $D \cap B_1 = \{p_{1,2} , p_{1,3}\}$ for some points $p_{1,2} \in B_1 \cap B_2$ and $p_{1,3} \in B_1 \cap B_3$.
    With $p_{1,2}$ and $p_{1,3}$ fixed, there is at most one curve in $\operatorname{Hyp}_{C_1}(B, 2)$ that contains both $p_{1,2}$ and $p_{1,3}$. Indeed, assume there are two such curves $D_{1}, D_{2} \in \operatorname{Hyp}_{C_{1}}(B,2)$, then we have 
    $
    D_{1} \cap B_{2} = \{p_{1,2}\} = D_{2} \cap B_{2}
    $
    and
    $
    D_{1} \cap B_{3} = \{p_{1,3}\} = D_{2} \cap B_{3}
    $.
    Then,
    $$
    (D_{1} \cdot B_{2})_{p_{1,2}} = (D_{1} \cdot B_{2}) = 2 \alpha_{2} + \beta_{2} = (D_{2} \cdot B_{2}) = (D_{2} \cdot B_{2})_{p_{1,2}}.
    $$
    
    Since $(D_{1} \cdot D_{2}) = C_{1}^{2} = 2$ and $D_{1} \cap D_{2} = \{p_{1,2}, p_{1,3}\}$, we must have
    $
    (D_{1} \cdot D_{2})_{p_{1,2}} = 1 = (D_{1} \cdot D_{2})_{p_{1,3}}
    $.
    This leads to
    $
    (D_{1} \cdot D_{2})_{p_{1,2}} 
    = 1 
    < 2 
    \leq 2 \alpha_{2} + \beta_{2} 
    = (D_{1} \cdot B_{2})_{p_{1,2}} 
    = (D_{2} \cdot B_{2})_{p_{1,2}}
    $.
    Note that $p_{1,2}$ is a smooth point on all three curves $D_{1}, D_{2}$, and $B_{2}$. But then, by Thm.~\ref{thm: Strong Triangle Inequality}, the smallest two among
    $$
    (D_{1} \cdot D_{2})_{p_{1,2}},
    \quad
    (D_{1} \cdot B_{2})_{p_{1,2}},
    \quad
    (D_{2} \cdot B_{2})_{p_{1,2}}
    $$
    are equal, contradiction.
    
    Since there are $(B_{1} \cdot B_{2})$ choices for $p_{1,2}$ and $(B_{1} \cdot B_{3})$ choices for $p_{1,3}$, we conclude that 
    $$
    \lvert \operatorname{Hyp}_{C_1}(B,2) \rvert 
    \leq (B_{1} \cdot B_{2}) (B_{1} \cdot B_{3})
    = (2\alpha_2 + \beta_2) (2\alpha_3 + \beta_3)
    $$
    in this case, and hence
    $$
    \lvert \mathcal{E}(B) \rvert 
    \leq 
    1 + (2 \alpha_{2} + \beta_{2})(2 \alpha_{3} + \beta_{3}) + \lvert N \rvert.
    $$

    If $\{B_1 , B_2 , B_3\} \cap \lvert C_1 \rvert = \{B_1 , B_2\}$, then we may switch $B_1$ and $B_2$. By symmetry, we obtain
    $$
    \lvert \operatorname{Hyp}_{C_1}(B,2) \rvert 
    \leq (B_{1} \cdot B_{2}) (B_{1} \cdot B_{3}) + (B_2 \cdot B_1) (B_2 \cdot B_3)
    = 4(2\alpha_3 + \beta_3),
    $$
    and hence
    $$
    \lvert \mathcal{E}(B) \rvert 
    \leq 1 + 4(2 \alpha_{3} + \beta_{3}) + \lvert N \rvert.
    $$

    If $\{B_1 , B_2 , B_3\} \subset \lvert C_1 \rvert$, then we may switch $B_2$ and $B_3$ with $B_1$. By symmetry, we obtain
    $$
    \lvert \operatorname{Hyp}_{C_1}(B,2) \rvert 
    \leq (B_{1} \cdot B_{2}) (B_{1} \cdot B_{3}) + (B_2 \cdot B_1) (B_2 \cdot B_3) + (B_3 \cdot B_1) (B_3 \cdot B_2)
    = 12.
    $$
    Note that $\lvert N \rvert = 6$ in this case, hence we obtain
    $$
    \lvert \mathcal{E}(B) \rvert 
    \leq 1 + 12 + \lvert N \rvert
    = 1 + 12 + 6
    = 19.
    $$
    This completes the proof.
\end{proof}

\begin{proposition}\label{prop: emptyness of exceptional set Fe}
    Consider $\mathbb{F}_e, e\geq 2$. Let $B$ be as in Prop.~\ref{prop: Fe Hyp(B,2)=E(B)}. 
    
    If $\alpha_i \geq 3$ for all $i \in \{1,2,3\}$, $\beta_l \geq 3$ for some $l \in \{1,2,3\}$, and $B$ is general, then
    $
    \mathcal{E}(B) = \emptyset
    $.
\end{proposition}
\begin{proof}
    By Prop.~\ref{prop: Fe Hyp(B,2)=E(B)}, we know $\operatorname{Hyp}(B,2) = \operatorname{Hyp}_{C_0} (B,2) \cup \operatorname{Hyp}_{f} (B,2)$.
    
    We first note that $f\vert_{B_i}$ is a base point free linear system on $B_i$ for any $i \in \{1,2,3\}$. For any $p \in B_i$, let 
    $f_p$ be the unique curve in $ \lvert f \rvert$ that contains $p$. If $B_i$ is smooth, then by Bertini's theorem (\cite[Ch.~\Romannum{3}, Cor.~10.9, p.~274]{Hartshorne1977}), we know for a general point $p \in B_i$ that $f_p$ intersects $B_i$ transversally. 
    Define the following set of points:
    $$
    T_f(B_i)
    \coloneqq
    \{
    p \in B_i
    \
    \vert
    \ 
    f_p
    \text{ is tangent to } B_i
    \}.
    $$
    Then $T_f(B_i)$ is a finite set.

    We choose $B$ general in the following sense: 
    $B_1 \in \lvert \alpha_1 C_1 + \beta_1 f \rvert$ is smooth; 
    $B_2 \in \lvert \alpha_2 C_1 + \beta_2 f \rvert$ is smooth, intersects $B_1$ transversally, and 
    $B_2 \cap T_f(B_1) = \emptyset$; 
    and $B_3 \in \lvert \alpha_3 C_1 + \beta_3 f \rvert$ is smooth, intersects $B_1 \cup B_2$ transversally,
    and 
    $B_3 \cap T_f(B_1) = \emptyset = B_3 \cap T_f(B_2)$. 
    Moreover, recall that we have assumed $\beta_l \geq 3$ for some $l\in \{1,2,3\}$, and we also require that $B_l$ intersects $C_0$ transversally.

    Then, by construction, 
    for any $i < j$, $i,j \in \{1,2,3\}$ and for any point $p_{i,j} \in B_i \cap B_j$, $f_{p_{i,j}}$ intersects $B_i$ transversally. We obtain
    $$
    \# f_{p_{i,j}} \cap B
    \geq 
    \# f_{p_{i,j}} \cap B_i
    = 
    (f_{p_{i,j}} \cdot B_i)
    = \alpha_i
    \geq 3.
    $$
    Thus $f_{p_{i,j}} \not\in \operatorname{Hyp}_{f} (B,2)$. By the first part of the proof of Thm.~\ref{thm: Fe bound E(B)}, we know 
    $
    \operatorname{Hyp}_{f} (B,2) = \emptyset.
    $
    Moreover, we have
    $$
    \# C_0 \cap B
    \geq 
    \# C_0 \cap B_l
    = 
    (C_0 \cdot B_l)
    = \beta_l
    \geq 3.
    $$
    Thus, $\operatorname{Hyp}_{C_0} (B,2) = \emptyset$.
    
    Combining with Prop.~\ref{prop: Fe Hyp(B,2)=E(B)}, we have
    $$
    \mathcal{E}(B) 
    = \operatorname{Hyp}(B,2)
    = \operatorname{Hyp}_{C_0} (B,2) \cup \operatorname{Hyp}_{f} (B,2)
    = \emptyset.
    $$
\end{proof}

\begin{remark} \label{rmk: without Corvaja--Zannier}
    In \cite[Thm.~1, Thm.~2]{CZ13}, Corvaja and Zannier bound the degree of a curve on a surface of log-general type which is also a finite affine smooth covering of $\mathbb{G}_{m}^{2}$ in terms of the number of intersection points with the boundary divisor and the geometric genus of the curve. We note that, in the case $S=\mathbb{F}_{e},e\geq 2$, if, in addition to the assumption in Prop.~\ref{prop: Fe Hyp(B,2)=E(B)}, none of the three components is (very) ample, then the Corvaja--Zannier method does not apply. Indeed, in this case the rational curve with negative self-intersection $C_{0}$ is contained in $\mathbb{F}_{e} \setminus B$, hence $\mathbb{F}_{e} \setminus B$ is not an affine surface. Of course, $C_0$ cannot be contracted to a point, for it would produce a surface singularity.

    From a birational geometric point of view, the open surfaces $\mathbb{G}_m^2$ and $\mathbb{P}^2\setminus B$, where $B$ is a three-component simple normal crossing divisor, all have log-irregularity $2$. However, if $B = B_1 \cup B_2 \cup B_3$ is a $3C-curve$ on a Hirzebruch surface $\mathbb{F}_e, e\geq 0$, such that $\mathcal{L}(B_1), \mathcal{L}(B_2), \mathcal{L}(B_3)$ generate $N^1(\mathbb{F}_e)_{\mathbb{Q}} \coloneqq \operatorname{Pic}(\mathbb{F}_e) \otimes_{\mathbb{Z}} \mathbb{Q}$ as a $\mathbb{Q}-$vector space (e.g., take $B_1 , B_2 \in \lvert C_1 \rvert$ and $B_3 \in \lvert C_1 + (e+1)f \rvert$), then the open surface $\mathbb{F}_e \setminus B$ has log-irregularity $1$. Indeed, consider the residue exact sequence
    $$
    0 \to \Omega_{\mathbb{F}_e} \to \Omega_{\mathbb{F}_e}(\log B) \to \oplus_{i=1}^3 \mathcal{O}_{B_i} \to 0.
    $$
    We have $h^0(\mathbb{F_e} , \Omega_{\mathbb{F}_e}) = 0$, and the image of the cycle class map $\mathbb{C}^{\oplus3} \cong \oplus_{i=1}^3H^0(B_i, \mathcal{O}_{B_i}) \to H^{1,1}(\mathbb{F}_e)$ has dimension $2$ by assumption. It follows that the log-irregularity of $\mathbb{F}_e\setminus B$ is $h^0 \left( \mathbb{F}_e , \Omega_{\mathbb{F}_e}(\log B) \right) = 1$. So the cases we study are not covered by \cite{caporaso2024hypertangencyplanecurvesalgebraic}, \cite{CZ08}, and \cite{CZ13}.
\end{remark}

\section{Hyper-bitangent Curves on \texorpdfstring{$\FF_{1}$}{F1}}
\label{sec:Hyper-bitangent Curves on F1}
In this section, we investigate $\operatorname{Hyp}(B,2)$ for the non-minimal Hirzebruch surface $S = \mathbb{F}_{1} \cong \operatorname{Bl}_{pt} \mathbb{P}^{2}$. It turns out to be a more subtle case than that of the minimal Hirzebruch surfaces.

Now, we consider a $3C$-curve $B = B_1 \cup B_2 \cup B_3$ on $\mathbb{F}_{1}$ such that $K_{\mathbb{F}_{1}} + \mathcal{L}(B)$ is big. Assume $B \in \lvert \alpha C_{1} + \beta f \rvert$ and $B_i \in \lvert \alpha_i C_1 + \beta_i f \rvert$ for $i=1,2,3$. By Cor.~\ref{cor: Irreducible Curves on Hirzebruch Surfaces}, if, say, $\beta_1 < 0$, then $\beta_1 = -1$, $\alpha_1 = 1$, and $B_1 = C_0$. As $B$ is reduced, it contains at most one copy of $C_0$. Hence $\beta_2 , \beta_3 \geq 0$. It follows that $\beta \geq -1$.

Since $K_{\mathbb{F}_{1}} = -2 C_{1} - f$, by Prop.~\ref{prop: Criterion for bigness on Hirzebruch Surfaces}, we have that $K_{\mathbb{F}_{1}} + \mathcal{L}(B)$ is big if and only if $\alpha > 2$ and $\alpha + \beta > 3$.

Assume that none of the irreducible components of $B$ belongs to $\lvert f \rvert$ or is $C_{0}$; then $\beta \geq 0$. Let us first consider the case $\beta = 0$. Then, $B \in \lvert \alpha C_{1} \rvert$ with $\alpha \geq 4$, and it does not contain any integral curve of type $u C_{1} + v f$ with $u \geq 0$ and $v > 0$ as a component. Hence, $B_{i} \in \lvert \alpha_{i} C_{1} \rvert$ with $\alpha_{i} \geq 1$ for all $i \in \{1, 2, 3\}$. Consider the blow-up morphism $\pi\colon \mathbb{F}_{1} \to \mathbb{P}^{2}$. We have $B = \pi^{*} \mathcal{O}(\alpha)$ with $\alpha \geq 4$. In this case, the exceptional curve is $E = C_{0}$ and $B \cap E = \emptyset$ since $(B \cdot E) = \sum_{i}(\alpha_{i} C_{1} \cdot C_{0}) = 0$. Thus, for any $D \in \operatorname{Hyp}(B,2)$, we have $\pi(D) \in \operatorname{Hyp}(\mathbb{P}^2,\pi(B),2) = \mathcal{E}(\mathbb{P}^2,\pi(B))$, where the last equality holds by \cite{caporaso2024hypertangencyplanecurvesalgebraic}. Therefore, $D$ must be a rational curve itself, and we see that the study of $\operatorname{Hyp}(B,2)$ in this case can be completely reduced to the case on $\mathbb{P}^{2}$, which has been studied in \cite{caporaso2024hypertangencyplanecurvesalgebraic}.

Recall that 
$N = \bigcup_{i < j} B_{i} \cap B_{j}$, and 
$\lvert N \rvert 
= \sum_{i < j} (B_i \cdot B_j)$ by the definition of a $3C$-curve. 
For the remaining cases with $\beta > 0$, we have the following proposition and theorems:
\begin{proposition}\label{prop: F1 Hyp-C1-(B,2)}
    Consider $\mathbb{F}_{1}$. Let $B = B_{1} \cup B_{2} \cup B_{3} \in \lvert \alpha C_{1} + \beta f \rvert$ be a $3C$-curve on $\mathbb{F}_{1}$, and denote $B_i \in \lvert \alpha_i C_1 + \beta_i f \rvert$ for $i=1,2,3$. Assume that $\beta>0$, $K_{\mathbb{F}_{1}} + \mathcal{L}(B)$ is big, and that none of $B_{1}$, $B_{2}$, $B_{3}$ is contained in $\lvert f \rvert$ or is $C_{0}$. Assume that $\beta_{3} \geq 1$ and 
    $\alpha_1 + \beta_1 
    \leq
    \alpha_2 + \beta_2$.

    Then,
    $$
    \lvert \operatorname{Hyp}_{C_1}(B,2) \rvert 
    \leq \frac{1}{2} \lvert N \rvert (\lvert N \rvert - 1).
    $$
\end{proposition}
\begin{proof}
    By assumption, we have $\alpha_{i} \geq 1$ and $\beta_{i} \geq 0$ for all $i = 1, 2, 3$, and $\beta_3 \geq 1$.
    
    Consider the blow-up morphism $\pi \colon \mathbb{F}_1 \to \mathbb{P}^2$, we have 
    $$
    h^0 (\mathbb{F}_1, C_1) 
    = h^0 \left( \mathbb{F}_1,\pi^{*}\mathcal{O}(1) \right)
    = h^0 \left(\mathbb{P}^2, \mathcal{O}(1) \right) 
    = 3.
    $$
    Therefore, through any two points on $\mathbb{F}_1$, there is at most one integral curve in $\lvert C_1 \rvert$ which contains both of them. Note that such a curve is a smooth rational curve.

    Let $D \in \operatorname{Hyp}_{C_1}(B,2)$, then $(D \cdot B_i) = \alpha_i + \beta_i \geq 1 > 0$ for any $i = 1,2,3$, which implies that $D \cap B_i \neq \emptyset$. As $B_1 \cap B_2 \cap B_3 = \emptyset$, we conclude that $\# (D \cap B) = 2$ and at least one of the two points in $D\cap B$ must be contained in $N$. Let $D\cap B = \{ p, q \}$ and $p \in N$. Assume $p \in B_{i_1} \cap B_{i_2}$, $i_1\neq i_2$; then $q \in B_{i_3}$, $\{i_1 , i_2 , i_3\} = \{1 , 2 , 3\}$. Then $D$ must intersect one of $B_{i_1}$ and $B_{i_2}$ transversally at $p$. Without loss of generality, let it be $B_{i_1}$.

    Write $\operatorname{Hyp}_{C_1}(B,2) = H_1 \cup H_2$, where
    $
    H_1 \coloneqq \left\{
    D \in \operatorname{Hyp}_{C_1}(B,2)
    \ \big\vert \ 
    D \cap B \subseteq N
    \right\}
    $
    and
    $
    H_2 \coloneqq \left\{
    D \in \operatorname{Hyp}_{C_1}(B,2)
    \ \big\vert \ 
    D \cap B \not\subseteq N
    \right\}
    $.

    If $H_2 = \emptyset$, then $D \in H_1$ and $D\cap B = \{p,q\} \subseteq N$. Since there is at most one integral curve in $\lvert C_1 \rvert$ which contains both $p$ and $q$, $D$ must be this curve. We have at most
    $
    \binom{\lvert N \rvert}{2} = \frac{1}{2} \lvert N \rvert (\lvert N \rvert - 1)
    $
    such curves, therefore the statement of the proposition holds in this case.

    It remains to study the case $H_2 \neq \emptyset$. 
    
    In this case, for $D \in H_2$, we have $D \cap B_{i_3} = \{ q \}$ and $D \cap B_{i_1} = D \cap B_{i_2} = \{ p \}$. Note that $p$ is a smooth point on both $B_{i_1}$ and $B_{i_2}$ by the definition of a $3C$-curve. By the assumption that $D$ is transverse to $B_{i_1}$ at $p$, we have
    $$
    1 
    = \operatorname{mult}_{p}(D)\operatorname{mult}_{p}(B_{i_1})
    = (B_{i_1} \cdot D)_{p}
    = (B_{i_1} \cdot D)
    = \alpha_{i_1} + \beta_{i_1},
    $$
    which implies that $\alpha_{i_1} = 1$ and $\beta_{i_1} = 0$. Since we have assumed 
    $\beta_3 \geq 1$ 
    and 
    $\alpha_1 + \beta_1 \leq \alpha_2 + \beta_2$, 
    we have $i_1 \in \{1,2\}$ and $B_1 \in \lvert C_1 \rvert$.

    By contracting $C_0$, we obtain the blow-up morphism $\pi \colon \mathbb{F}_1 \to \mathbb{P}^2$, and
    $$
    \pi(D) , \pi(B_{i_1}) \in \lvert \mathcal{O}(1) \rvert,
    \quad
    \pi(B_{i_2}) \in \lvert \mathcal{O}(\alpha_{i_2} + \beta_{i_2}) \rvert, 
    \quad 
    \pi(B_{i_3}) \in \lvert \mathcal{O}(\alpha_{i_3} + \beta_{i_3}) \rvert.
    $$
    Since $(C_1 \cdot C_0) = 0$, we know $B_{i_1} \cap C_0 = \emptyset = D \cap C_0$ for any $D \in \operatorname{Hyp}_{C_1}(B,2)$. Hence, $\pi(B)$ has normal crossing singularity at $\pi(p)$, and for any $D \in \operatorname{Hyp}_{C_1}(B,2)$, $\pi(D) \subset \mathbb{P}^2$ is a line which does not contain the blow-up point.

    Then, we have the following two subcases:
    \begin{enumerate}[label=(\alph*), ref=(\alph*)]
        \item $\{B_1 , B_2\} \subset \lvert C_1 \rvert$.
        
        In this subcase, $\lvert N \rvert = 1 + 2(\alpha_3 + \beta_3)$. For $D \in H_1$, it cannot coincide with $B_1$ or $B_2$, hence one of $p$ and $q$ must be contained in $B_1 \cap B_3$, and the other one must be contained in $B_2 \cap B_3$. Therefore,
        $$
        \lvert H_1 \rvert \leq (B_1\cdot B_3) (B_2 \cdot B_3) = (\alpha_3 + \beta_3)^2.
        $$
        
        Now consider $D \in H_2$. As discussed before, we have $D \cap B_{i_1} = D \cap B_{i_2} = \{ p \}$, $D \cap B_{i_3} = \{ q \}$, $q\not\in N$, and $B_{i_1} \in \lvert C_1\rvert$ intersects $D$ transversally at $p$.
            
        If $p \in B_1 \cap B_3$, then $i_1 = 1$ and $D \cap B_3 = \{p\}$. Hence $\pi(D)$ must be the tangent line to $\pi(B_3)$ at $\pi(p)$, and there are at most $(B_1 \cdot B_3) = \alpha_3 + \beta_3$ such lines. 
            
        If $p \in B_2 \cap B_3$, then $i_1 = 2$ and $D \cap B_3 = \{p\}$. By symmetry, there are at most $\alpha_3 + \beta_3$ choices for such $D$.
            
        If $p \in B_1 \cap B_2$, then $\{ \pi(p) \} = \pi(B_1) \cap \pi(B_2)$, and $\pi(D)$ is a line containing $\pi(p)$ which is hypertangent to $\pi(B_3)$ at $\pi(q)$. The lines containing $\pi(p)$ give a morphism of degree $\alpha_3 + \beta_3$ from $\pi(B_3)$ to $\mathbb{P}^1$, and those lines which are hypertangent to $\pi(B_3)$ correspond to the points of maximal ramification index $\alpha_3 + \beta_3  - 1$ of this morphism. Denote $b \coloneqq \alpha_3 + \beta_3$; by the Riemann-Hurwitz theorem for singular curves (see, for instance, \cite[p.~2]{GarciaLax96}), we obtain
        \begin{align}
            \#\text{points of maximal ramification index}
            & \leq 
            \frac{2p_a\left( \pi(B_3) \right) - 2 - b \left( 2p_{a}(\mathbb{P}^1) - 2 \right)}{b- 1} \\
            & = \frac{(b - 1)(b - 2) - 2 + 2b }{b - 1} = b. 
        \end{align}
        Hence there are at most $b = \alpha_3+\beta_3$ such lines.
        
        Therefore, we obtain
        $$
        \lvert H_2 \rvert 
        \leq 3 (\alpha_3 + \beta_3). 
        $$
            
        Summing up, we obtain
        \begin{align}
            \lvert \operatorname{Hyp}_{C_1}(B,2) \rvert
            & \leq (\alpha_3 + \beta_3)^2 + 3 (\alpha_3 + \beta_3) \\
            & \leq 2(\alpha_3 + \beta_3)^2 + (\alpha_3 + \beta_3) = \frac{1}{2} \lvert N \rvert (\lvert N \rvert - 1)
        \end{align}
        in the subcase $\{B_1,B_2\} \subset \lvert C_1 \rvert$ as desired.

        \item $\{B_1 , B_2\} \cap \lvert C_1 \rvert = \{B_1\}$.
        
        In this subcase, we must have $i_1 = 1$. Denote $a \coloneqq \alpha_2 + \beta_2 \geq 2$ and $b \coloneqq \alpha_3 + \beta_3 \geq 2$; then
        $$
        (B_1 \cdot B_2) = a,
        \quad
        (B_1 \cdot B_3) = b,
        $$
        and
        $$
        (B_2 \cdot B_3) 
        = \alpha_2\alpha_3 + \alpha_2\beta_3 + \alpha_3\beta_2
        = \alpha_2 b + \alpha_3\beta_2
        = ab - \beta_2\beta_3.
        $$
        Hence, we obtain 
        $\lvert N \rvert 
        = a + b + ab - \beta_2\beta_3$. 
        
        For $D \in H_1$, it cannot contain points in $B_1 \cap B_2$ and $B_1\cap B_3$ at the same time; otherwise, it coincides with $B_1$. Therefore, one of $p$ and $q$ must be contained in $B_2 \cap B_3$, and the other one must be contained either in $B_1 \cap B_2$ or $B_1 \cap B_3$. We obtain 
        $$
        \lvert H_1 \rvert 
        \leq \left[(B_1 \cdot B_2) + (B_1 \cdot B_3)\right] (B_2 \cdot B_3)
        \leq (a+b)(ab - \beta_2\beta_3).
        $$
    
        For $D \in H_2$, we have $\{ p \} = B_1 \cap D = B_{i_2} \cap D$. Since $i_2 \in \{ 2 , 3 \}$, $\alpha_2 + \beta_2 \geq 2$ and $\alpha_3 + \beta_3 \geq 2$, $\pi(D)$ must be the tangent line to $\pi(B_{i_2})$ at $\pi(p)$. There are at most $(B_1 \cdot B_2) + (B_1 \cdot B_3)$ such lines, hence
        $$
        \lvert H_2 \rvert \leq (B_1 \cdot B_2) + (B_1 \cdot B_3) = a + b.
        $$
    
        Summing up, we obtain
        $$
        \lvert \operatorname{Hyp}_{C_1}(B,2) \rvert 
        = \lvert H_1 \rvert + \lvert H_2 \rvert
        \leq (a+b)(ab - \beta_2\beta_3 +1)
        = (a+b)(\lvert N \rvert + 1 -a-b).
        $$
        \begin{claim}
            $
            (a+b)(\lvert N \rvert + 1 -a-b)
            <
            \frac{1}{2}\lvert N \rvert(\lvert N \rvert - 1).
            $
        \end{claim}
        \begin{proof}(of the claim)
        Denote $c \coloneqq a+b$. Then $c \geq 4$, and 
        $$
        \lvert N \rvert 
        = a + b + ab - \beta_2\beta_3
        = c + \alpha_2 b +\alpha_3\beta_2
        \geq c + b
        \geq c + 2.
        $$
        The statement in the claim is equivalent to
        $$
        \lvert N \rvert^2 
        - \left( 2c+1\right) \lvert N \rvert
        + 2c(c-1)
        > 0.
        $$
        We have
        \begin{align*}
            \lvert N \rvert^2 
            - \left( 2c+1\right) \lvert N \rvert
            + 2c(c-1)
            & = \lvert N \rvert^2 
            - \left( 2c+1\right) \lvert N \rvert
            + c(c+1) + c(c-3) \\
            & = \left( \lvert N \rvert - c \right) \left( \lvert N\rvert - c - 1 \right) + c(c-3)
            > 0
        \end{align*}
        as $N \geq c+2$ and $c \geq 4$. Hence the statement in the claim holds.
        \end{proof}
        From the claim, it immediately follows that
        $$
        \lvert \operatorname{Hyp}_{C_1}(B,2) \rvert
        < \frac{1}{2} \lvert N \rvert (\lvert N \rvert - 1)
        $$
        in the subcase 
        $\{B_1 , B_2\} \cap \lvert C_1 \rvert
        = \{B_1\}$.    
    \end{enumerate}
    
    This completes the proof.
\end{proof}

\begin{proposition}\label{prop: F1 Hyp(B,2)=E(B)}
    Consider $\mathbb{F}_1$. Under the same assumptions as in Prop.~\ref{prop: F1 Hyp-C1-(B,2)}, consider $d_1 C_1 + d_2 f \in \operatorname{Pic}(\mathbb{F}_1)$ with $d_1 \geq 1$, $d_2 \geq 0$ and $d_1 + d_2 \geq 2$.

    If $\operatorname{Hyp}_{d_{1}C_{1}+d_{2}f}(B,2) \neq \emptyset$, then the following holds:

    \begin{enumerate}
        \item\label{prop: F1 Hyp(B,2)=E(B) (1)} 
        $B_1 \in \lvert C_1 \rvert$.
        
        \item\label{prop: F1 Hyp(B,2)=E(B) (2)} For any $d_1 \geq 2$, let $\mathfrak{I} \subseteq \{1,2\}$ be the set of indices such that $B_i \in \lvert C_1 \rvert$ for all $i \in \mathfrak{I}$. Then
        $$
        \operatorname{Hyp}_{d_1 C_1}(B,2) = \bigcup\limits_{\substack{
        i \in \mathfrak{I} \\
        \{i,j\} = \{1,2\} \\
        p_{1,2} \in B_1 \cap B_2 \\ p_{i,3}\in B_i \cap B_3}} \operatorname{Hyp}_{d_{1} C_{1}}^{d_{1} - 1}(B_{j}, p_{1,2}) \cap \operatorname{Hyp}_{d_{1} C_{1}}(B_{3}, p_{i,3}),
        $$
        and for any $d_1 \geq 3$, $\operatorname{Hyp}_{d_1 C_1}(B,2) \neq \emptyset$ 
        only if 
        $\{B_1 , B_2\} \subset \lvert C_1 \rvert$.
        Moreover, $\left\lvert \bigcup_{d_1 \geq 2} \operatorname{Hyp}_{d_1 C_{1}}(B,2) \right\rvert$ is finite and is effectively bounded in terms of $\alpha_i$'s and $\beta_i$'s.
        \item\label{prop: F1 Hyp(B,2)=E(B) (3)} If $d_2 > 0$, then $d_2 = 1$. Furthermore, let $\mathfrak{I} \subseteq \{1,2\}$ be the set of indices such that $B_i \in \lvert C_1 \rvert$ for all $i \in \mathfrak{I}$. Then
        $$
        \operatorname{Hyp}_{d_1 C_1 + f}(B,2)
        = \bigcup\limits_{\substack{
        i \in \mathfrak{I} \\
        \{i,j\} = \{1,2\} \\ 
        p_{1,2} \in B_1 \cap B_2 \\ 
        p_{i,3} \in B_i \cap B_3} } \operatorname{Hyp}_{d_{1} C_{1} + f}^{d_{1}}(B_{j}, p_{1,2}) \cap \operatorname{Hyp}_{d_1 C_1 + f}(B_{3}, p_{i,3}).
        $$
        Moreover, for any $d_1 \geq 2$,
        $\operatorname{Hyp}_{d_1 C_1 + f}(B,2) \neq \emptyset$
        only if
        $\{B_1 , B_2\} \subset \lvert C_1 \rvert$. 
    \end{enumerate}
    In particular, $\operatorname{Hyp}(B,2) = \mathcal{E}(B)$.
\end{proposition}

\begin{proof}
    By assumption, we have $\alpha_{i} \geq 1$ and $\beta_{i} \geq 0$ for all $i = 1, 2, 3$, and $\beta_3 \geq 1$.
    
    Let $D$ be a curve in $\operatorname{Hyp}_{d_1 C_1 + d_2 f}(B,2)$, where $d_1 \geq 1$, $d_2 \geq 0$ and $d_1 + d_2 \geq 2$. Then $D \cap B_i \neq \emptyset$ for all $i = 1,2,3$. Since $B_{1} \cap B_{2} \cap B_{3} = \emptyset$, it follows that $\#(D \cap B) = 2$. This implies that any point in $D\cap B$ is a unibranch point on $D$.

    \begin{claim}
        $D \cap B_{3} \subseteq N$ and $\#(D \cap B_{3}) = 1$.
    \end{claim}
    \begin{proof}[Proof of the Claim]
        Suppose there exists a point $q$ in $D \cap B_{3}$ that is not contained in $N$. Then, $D$ must intersect $B_{1} \cup B_{2}$ in exactly one point. Since 
        $(B_i \cdot B_j) 
        = \alpha_i \alpha_j + \alpha_i \beta_j + \alpha_j \beta_i 
        \geq 1 > 0$, 
        we know that $B_{i} \cap B_{j} \neq \emptyset$ for any $i \neq j$. Thus, $D$ intersects $B_1 \cup B_{2}$ at one point $p_{1,2} \in B_1 \cap B_2$, which must be a unibranch $n$-fold point of $D$ for some $n \geq 1$.

        Then, $D$ is transverse to one of $B_{1}$ and $B_{2}$ at $p_{1,2}$. Assume $D$ is transverse to $B_{i}$ at $p_{1,2}$, $i \in \{1, 2\}$; then
        $$
        n 
        = (D \cdot B_{i}) 
        = \left( (d_{1} C_{1} + d_2 f) \cdot (\alpha_{i} C_{1} + \beta_{i} f) \right) 
        = (\alpha_{i} + \beta_i) d_{1} + \alpha_{i} d_2 
        \geq d_1 + d_2 
        \geq d_1.
        $$
        By Lem.~\ref{lem: Upper bound of unibranch multiplicity on Hirzebruch Surfaces}, we see all the equalities here must hold, which gives $\alpha_i = 1$, $\beta_i = 0$ and $d_2 = 0$. But then $D \in \lvert d_1 C_1 \rvert$, and 
        $n = d_1 = d_1 + d_2 \geq 2$
        contradicts Lem.~\ref{lem: Special case of bound of unibranch multiplicity on F1}. Therefore, we have proved that $D \cap B_3 \subseteq N$.

        Now, if $\#(D \cap B_{3}) = 2$, we have $D \cap B_{3} = \{p_{1,3}, p_{2,3}\}$ for some points $p_{1,3} \in B_1 \cap B_3$ and $p_{2,3} \in B_2 \cap B_3$. Then, $D$ is transverse to $B_{i}$ or $B_{3}$ at $p_{i,3}$ for any $i = 1, 2$.

        If $D$ is transverse to $B_{i}$ at $p_{i,3}$ for some $i \in \{1,2\}$, as before, we may assume $p_{i,3}$ is a unibranch $n$-fold point of $D$ for some $n \geq 1$. Then,
        $$
        n 
        = (D \cdot B_{i})_{p_{i,3}} 
        = (D \cdot B_{i}) 
        = (\alpha_{i} + \beta_{i}) d_{1} + \alpha_i d_2 
        \geq d_{1} + d_2 
        \geq d_1.
        $$
        Again, by Lem.~\ref{lem: Upper bound of unibranch multiplicity on Hirzebruch Surfaces}, we see all the equalities here must hold, which implies $\alpha_i = 1$, $\beta_i = 0$ and $d_2 = 0$. But then $D \in \lvert d_1 C_1 \rvert$, and 
        $n = d_1 = d_1 + d_2 \geq 2$ 
        contradicts Lem.~\ref{lem: Special case of bound of unibranch multiplicity on F1}.

        Hence, $D$ is transverse to $B_3$ at both $p_{1,3}$ and $p_{2,3}$. Denote $m_{i} \coloneqq \operatorname{mult}_{p_{i,3}}(D)$ for any $i = 1, 2$. Then by Lem.~\ref{lem: Upper bound of unibranch multiplicity on Hirzebruch Surfaces}, we know $m_i \leq d_1$ for $i = 1,2$. We obtain
        \begin{multline}
        2d_1 \geq m_1 + m_2 
        = (D \cdot B_3)_{p_{1,3}} + (D \cdot B_3)_{p_{2,3}}
        \\
        = (D \cdot B_3) 
        = (\alpha_3 + \beta_3) d_1 + \alpha_3 d_2 
        \geq 2d_1 + d_2 
        \geq 2d_1.
        \end{multline}
        Therefore, all the equalities here must hold, which gives 
        $$d_1 = m_1 = m_2 ,
        \quad
        \alpha_3 = \beta_3 = 1 ,
        \quad
        \text{and}
        \quad
        d_2 = 0.
        $$ But then $D \in \lvert d_1 C_1 \rvert$, and $d_1 = d_1 + d_2 \geq 2$. By Lem.~\ref{lem: Special case of bound of unibranch multiplicity on F1}, we must have $d_1 > m_1$ and $d_1 > m_2$, contradiction. This completes the proof of the claim.
    \end{proof}

    From the claim, we see that $D \cap B_3 = \{ p_{i,3} \}$ for some $i \in \{1, 2\}$. Notice that $D$ must be tangent to $B_3$ at $p_{i,3}$ because
    $
    (D \cdot B_{3}) 
    = (\alpha_{3} + \beta_{3}) d_{1} + \alpha_3 d_2 
    \geq 2 d_{1} + d_2 
    > d_{1} 
    \geq \operatorname{mult}_{p_{i,3}}(D)
    $.

    Thus, $D$ must be transverse to $B_{i}$ at $p_{i,3}$. Then, $D$ must meet $B_{i}$ at a further point. Indeed, we have
    $$
    (D \cdot B_{i}) 
    = (\alpha_{i} + \beta_{i}) d_{1} + \alpha_i d_2 
    \geq d_{1} + d_2 
    \geq d_1 
    \geq \operatorname{mult}_{p_{i,3}}(D) = (D \cdot B_{i})_{p_{i,3}}
    .$$
    If $D \cap B_i = \{ p_{i,3} \}$, then all the equalities here must hold, which gives $\alpha_i = 1$, $\beta_i = 0$, $d_2 = 0$ and 
    $\operatorname{mult}_{p_{i,3}}(D)
    = d_1 = d_1 + d_2 \geq 2$. 
    But then $D \in \lvert d_1 C_1 \rvert$, and by Lem.~\ref{lem: Special case of bound of unibranch multiplicity on F1} we know $\operatorname{mult}_{p_{i,3}}(D) < d_1$, contradiction.

    Since $D$ must meet the other component $B_{j}$ (where $j \neq i, 3$), we know there exists a point $p_{1,2} \in B_{1} \cap B_{2}$ such that $D \cap B_{j} = \{p_{1,2}\}$ and $D \cap B = \{p_{1,2}, p_{i,3}\}$.

    Now, $B_{j}$ and $D$ meet only at $p_{1,2}$. We know that $D$ cannot be transverse to $B_{j}$ at $p_{1,2}$ for the same reason as before. Thus, $D$ must be transverse to $B_{i}$ at $p_{1,2}$.

    Hence, $D \cap B_{i} = \{p_{1,2}, p_{i,3}\}$, and they intersect transversally at both of these points. By Lem.~\ref{lem: Upper bound of unibranch multiplicity on Hirzebruch Surfaces}, we know $d_1 \geq \operatorname{mult}_{p_{1,2}}(D)$ and $d_1 \geq \operatorname{mult}_{p_{i,3}}(D)$. We obtain
    \begin{multline}
        2d_{1}
        \geq \operatorname{mult}_{p_{1,2}}(D) + \operatorname{mult}_{p_{i,3}}(D)
        = (D \cdot B_{i})_{p_{1,2}} + (D \cdot B_{i})_{p_{i,3}} \\
        = (D \cdot B_{i})
        = (\alpha_{i} + \beta_{i}) d_{1} + \alpha_i d_2
        \geq d_1 + d_2.
    \end{multline}
    
    By assumption, we have $\alpha_{i} \geq 1$ and $\beta_{i} \geq 0$. Thus, we must have $\alpha_{i} = 1$, $\beta_{i} = 0$ and $d_2 \leq d_1$. Therefore, $B_{i}$ is contained in $\lvert C_{1}\rvert$. By the assumption 
    $\alpha_1 + \beta_1 \leq \alpha_2 + \beta_2$,
    we must have $B_1 \in \lvert C_1 \rvert$. This proves part (\ref{prop: F1 Hyp(B,2)=E(B) (1)}) of the proposition.
    
    Without loss of generality, we may assume $i = 1$. Then, $D\cap B_3 = \{p_{1,3}\}$ and $D \cap B_{1} = \{p_{1,2}, p_{1,3}\}$. Now, we have
    $$
    \operatorname{mult}_{p_{1,2}}(D) + \operatorname{mult}_{p_{1,3}}(D) 
    = (D \cdot B_{1})_{p_{1,2}} + (D \cdot B_{1})_{p_{1,3}} 
    = (D \cdot B_{1}) 
    = d_{1} + d_{2}.
    $$

    Let $m \coloneqq \operatorname{mult}_{p_{1,2}}(D)$. Then, $m \leq d_{1}$ by Lem.~\ref{lem: Upper bound of unibranch multiplicity on Hirzebruch Surfaces}, and we have
    $$
    D \in
    \operatorname{Hyp}_{d_{1} C_{1} + d_{2} f}^{m}(B_{2}, p_{1,2}) \cap \operatorname{Hyp}_{d_{1} C_{1} + d_{2} f}^{d_{1} + d_{2} - m}(B_{3}, p_{1,3}).
    $$

    \noindent \textbf{Case 1, $d_{2} = 0$:}
    In this case, we have $d_1 \geq 2$, $m < d_{1}$ by Lem.~\ref{lem: Special case of bound of unibranch multiplicity on F1}, and
    $$
    D \in \operatorname{Hyp}_{d_{1} C_{1}}^{m}(B_{2}, p_{1,2}) \cap \operatorname{Hyp}_{d_{1} C_{1}}^{d_{1} - m}(B_{3}, p_{1,3}).
    $$

    By Thm.~\ref{thm: bounddeltainvariant}, we have
    $$
    \delta_{D}(p_{1,2}) 
    \geq \frac{(m - 1) \left( (D \cdot B_{2}) - 1 \right)}{2} 
    = \frac{(m - 1) \left( \alpha_{2} d_{1} + \beta_{2} d_{1} - 1 \right)}{2}
    $$
    and
    $$
    \delta_{D}(p_{1,3}) 
    \geq \frac{(d_{1} - m - 1) \left( (D \cdot B_{3}) - 1 \right)}{2} 
    = \frac{(d_{1} - m - 1) \left( \alpha_{3} d_{1} + \beta_{3} d_{1} - 1 \right)}{2}.
    $$

    From these, we obtain
    \begin{equation}\label{ineq: F1 case1 genusineq pre}
        \begin{split}
            0 \leq p_{a}(D^{\nu}) 
            &\leq p_{a}(D) - \delta_{D}(p_{1,2}) - \delta_{D}(p_{1,3}) \\
            & \leq \frac{1}{2}(d_{1} - 1)(d_{1} - 2) - \frac{1}{2}(m - 1)(\alpha_{2} d_{1} + \beta_{2} d_{1} - 1) \\
            &\quad - \frac{1}{2}(d_{1} - m - 1)(\alpha_{3} d_{1} + \beta_{3} d_{1} - 1).
        \end{split}
    \end{equation}

    As $\alpha_{2} \geq 1$ and $\beta_{2} \geq 0$, we have
    \begin{equation}\label{ineq: F1 case1 genusineq}
        \begin{split}
            0 \leq p_{a}(D^{\nu}) 
            & \leq \frac{1}{2} (d_{1} - 1)(d_{1} - 2) - \frac{1}{2}(m - 1)(d_{1} - 1) 
            - \frac{1}{2}(d_{1} - m - 1)(\alpha_{3}d_{1} + \beta_{3}d_{1} - 1) \\
            & = \frac{1}{2}(d_{1} - 1)(d_{1} - m - 1)
            - \frac{1}{2}(d_{1} - m - 1)(\alpha_{3}d_{1} + \beta_{3}d_{1} - 1) \\
            & = -\frac{1}{2}(d_{1} - m - 1)(\alpha_{3} + \beta_{3} - 1)d_{1} 
            \leq 0.
        \end{split}
    \end{equation}

    Therefore, the equalities here must hold, which implies that $d_{1} = m + 1$ as claimed in part (\ref{prop: F1 Hyp(B,2)=E(B) (2)}) of the proposition.  Substituting this back into the inequality \eqref{ineq: F1 case1 genusineq pre}, we obtain
    \begin{align*}
        0 \leq p_{a}(D^{\nu}) 
        & \leq \frac{1}{2} m (m - 1) - \frac{1}{2}(m - 1)\left( (\alpha_{2} + \beta_{2})m + (\alpha_{2} + \beta_{2} - 1) \right) \\
        &= -\frac{1}{2} (m - 1)^{2} (\alpha_{2} + \beta_{2} - 1) \leq 0.
    \end{align*}
    Therefore, all the equalities must hold and $p_{a}(D^{\nu}) = 0$. 
    
    If $m \geq 2$ (equivalently, $d_1 \geq 3$), we must have $\alpha_{2} = 1$ and $\beta_{2} = 0$, that is, when $d_1 \geq 3$,
    $\operatorname{Hyp}_{d_1 C_1}(B,2) \neq \emptyset$
    only if $\{B_1 , B_2\} \subset \lvert C_1 \rvert$.
    Moreover, if 
    $\{B_1 , B_2\} \subset \lvert C_1 \rvert$,
    we also need to consider those curves in 
    $\operatorname{Hyp}_{d_1 C_1}(B,2)$ that intersect $B_2$ at two points. To do this, we may switch $B_1$ and $B_2$, then by symmetry, the preceding arguments remain valid.

    So, to prove the last statement of part (\ref{prop: F1 Hyp(B,2)=E(B) (2)}) of the proposition, we first need to study the cardinality of 
    $\operatorname{Hyp}_{2C_1}(B,2)$.

    \begin{claim}\label{claim: F1 Hyp 2C1 effectively bounded}
        With $p_{1,2}$ and $p_{1,3}$ fixed, there is at most one integral curve in $\operatorname{Hyp}_{2C_1}(B,2)$ that contains both $p_{1,2}$ and $p_{1,3}$
    \end{claim}
    \begin{proof}[Proof of the Claim]
        Let $D$ and $D^{\prime}$ be two such curves. Then as we have seen before, we have 
        $D \cap B_2 
        = \{ p_{1,2} \} 
        = D^{\prime} \cap B_2$ 
        and 
        $D \cap B_3 
        = \{ p_{1,3} \}
        = D^{\prime} \cap B_3$. 
        We obtain
        $$
        (D\cdot B_2)_{p_{1,2}} 
        = (D \cdot B_2) 
        = 2 (\alpha_2 + \beta_2)
        = (D^{\prime} \cdot B_2)
        = (D^{\prime} \cdot B_2)_{p_{1,2}}
        $$
        and
        $$
        (D\cdot B_3)_{p_{1,3}} 
        = (D \cdot B_3) 
        = d_1 (\alpha_3 + \beta_3)
        = (D^{\prime} \cdot B_3)
        = (D^{\prime} \cdot B_3)_{p_{1,3}}.
        $$
        
        Since $(D^{\prime} \cdot D) = 4$, we must have 
        $$
        (D^{\prime} \cdot D)_{p_{1,3}} 
        < 4
        \leq 2(\alpha_3 + \beta_3)
        = (D^{\prime} \cdot B_3)_{p_{1,3}}
        = (D \cdot B_3)_{p_{1,3}}.
        $$
        Note that $p_{1,3}$ is a smooth point on all three curves $D$, $D^{\prime}$ and $B_3$. But then, by Thm.~\ref{thm: Strong Triangle Inequality}, the smallest two among
        $$
        (D^{\prime} \cdot D)_{p_{1,3}}, 
        \quad
        (D^{\prime} \cdot B_3)_{p_{1,3}},
        \quad
        (D \cdot B_3)_{p_{1,3}}
        $$
        are equal, contradiction.
    
    \end{proof}
    
    If $\{B_1 , B_2\} \cap \lvert C_1 \rvert = \{B_1\}$, then since there are only $(B_1 \cdot B_2) = \alpha_2 + \beta_2$ choices for $p_{1,2}$ and $(B_1 \cdot B_3) = \alpha_3 + \beta_3$ choices for $p_{1,3}$, we obtain 
    $$\lvert \operatorname{Hyp}_{2C_1}(B,2) \rvert 
    \leq (\alpha_2 + \beta_2) (\alpha_3 + \beta_3).
    $$  
    
    If $\{B_1 , B_2\} \subset \lvert C_1 \rvert$, then 
    $B_1 \cap B_2 = \{p_{1,2}\}$.
    As mentioned before, we also need to consider those curves in 
    $\operatorname{Hyp}_{2 C_1}(B,2)$ that intersect $B_2$ at two points. By symmetry we obtain
    $$
    \lvert \operatorname{Hyp}_{2C_1}(B,2) \rvert 
    \leq (B_1 \cdot B_3) + (B_2 \cdot B_3) 
    = 2(\alpha_3 + \beta_3).
    $$ 
    
    In both cases, $\lvert \operatorname{Hyp}_{2C_1}(B,2) \rvert$ is effectively bounded.

    Finally, if $\{B_1 , B_2\} \cap \lvert C_1 \rvert = \{B_1\}$, we have seen that $\lvert \operatorname{Hyp}_{d_1 C_1} (B,2) \rvert = \emptyset$ for all $d_1 > 2$. If $\{B_1 , B_2\} \subset \lvert C_1 \rvert$, then by contracting $C_0$ we obtain the blow-up morphism $\pi\colon \mathbb{F}_{1} \to \mathbb{P}^{2}$. Note that $C_0 \cap B_1 = \emptyset = C_0 \cap B_2$. Therefore, the plane curve $\pi(B) = \pi(B_1) \cup \pi(B_2) \cup \pi(B_3)$ is still a $3C$-curve, $\pi(D) \in \operatorname{Hyp}(\mathbb{P}^2 , \pi(B) , 2)$, and $\pi(D) \cap \pi(B_{j})= \{ \pi(p_{1,j}) \}$ for $j = 2, 3$. We also have
    $$
    \pi(B_1),\pi(B_2) \in \lvert \mathcal{O}(1) \rvert,
    \quad
    \pi(B_3) \in \lvert \mathcal{O}(\alpha_3 + \beta_3) \rvert.
    $$
    Hence $\deg \pi(B) \geq 4$; we reduce the problem to the cases on $\mathbb{P}^2$, which are already studied in \cite[Prop.~3.2.1, Thm.~3.3.1, Thm.~3.3.2]{caporaso2024hypertangencyplanecurvesalgebraic}. From these, we conclude the proof of part (\ref{prop: F1 Hyp(B,2)=E(B) (2)}) of the proposition.

    \noindent \textbf{Case 2, $d_{2} > 0$:}
    Recall that we assumed that $B_1$ intersects $D$ transversally at $p_{1,2}$ and $p_{1,3}$, and we have
    $$
    \operatorname{mult}_{p_{1,2}}(D) + \operatorname{mult}_{p_{1,3}}(D) = (D \cdot B_{1})_{p_{1,2}} + (D \cdot B_{1})_{p_{1,3}} = (D \cdot B_{1}) = (D \cdot C_{1}) = d_{1} + d_{2},
    $$
    $m \coloneqq \operatorname{mult}_{p_{1,2}}(D) \leq d_1$, and
    $$
    D \in \operatorname{Hyp}_{d_{1} C_{1} + d_{2} f}(B,2) \subseteq \operatorname{Hyp}_{d_{1} C_{1} + d_{2} f}^{m}(B_{2}, p_{1,2}) \cap \operatorname{Hyp}_{d_{1} C_{1} + d_{2} f}^{d_{1} + d_{2} - m}(B_{3}, p_{1,3}).
    $$

    By Thm.~\ref{thm: bounddeltainvariant}, we have
    $$
    \delta_{D}(p_{1,2}) \geq \frac{(m - 1) \left( (D \cdot B_{2}) - 1 \right)}{2} = \frac{(m - 1) \left( \alpha_{2} d_{1} + \alpha_{2} d_{2} + \beta_{2} d_{1} - 1 \right)}{2}
    $$
    and
    $$
        \delta_{D}(p_{1,3}) 
        \geq \frac{(d_{1} + d_{2} - m - 1) \left( (D \cdot B_{3}) - 1 \right)}{2}
        = \frac{(d_{1} + d_{2} - m - 1) \left( \alpha_{3} d_{1} + \alpha_{3} d_{2} + \beta_{3} d_{1} - 1 \right)}{2}.    
    $$

    From these, we obtain
    \begin{equation}\label{ineq: F1 case2 genusineq pre}
        \begin{split}
            0 \leq p_{a}(D^{\nu}) & \leq p_{a}(D) - \delta_{D}(p_{1,2}) - \delta_{D}(p_{2,3}) \\
            & \leq \frac{1}{2}(d_{1}-1)(d_{1}+2d_{2}-2)
            - \frac{1}{2}(m-1) \left( \alpha_{2}d_{1}+\alpha_{2}d_{2}+\beta_{2}d_{1}-1\right) \\
            & \quad - \frac{1}{2}(d_{1}+d_{2}-m-1) \left( \alpha_{3}d_{1}+\alpha_{3}d_{2}+\beta_{3}d_{1}-1\right) 
            .
        \end{split}
    \end{equation}

    As $\alpha_{2}\geq 1, \beta_{2}\geq 0$, we get
        \begin{equation}\label{ineq: F1 case2 genusineq}
            \begin{split}
                0  \leq p_{a}(D^{\nu}) & \leq \frac{1}{2}(d_{1}-1)(d_{1}+2d_{2}-2) - \frac{1}{2}(m-1)(d_{1}+d_{2}-1) \\
                & \quad - \frac{1}{2} (d_{1}+d_{2}-m-1)( (\alpha_{3}+\beta_{3})d_{1} + \alpha_{3}d_{2} - 1 )  \\
                & = \frac{1}{2}
                \left( -(\alpha_{3}+\beta_{3}-1)d_{1}^{2} - (2\alpha_{3}+\beta_{3}-2)d_{1}d_{2} - \alpha_{3}d_{2}^{2}\right) \\
                & \quad + \frac{1}{2}\left( (\alpha_{3}+\beta_{3}-1)(m+1)d_{1} + ((\alpha_{3}-1)m+\alpha_{3})d_{2} \right)\\
                & = - \frac{1}{2}\left[ (\alpha_{3}+\beta_{3}-1)d_{1} + (\alpha_{3}-1)d_{2} \right](d_{1}-m) \\
                & \quad - \frac{1}{2} \left[(\alpha_{3}+\beta_{3}-1)d_{1} + \alpha_{3}d_{2} \right] (d_{2}-1)
                .
            \end{split}
        \end{equation}

    Note that $d_{1} \geq m$ by Lem.~\ref{lem: Upper bound of unibranch multiplicity on Hirzebruch Surfaces} and $d_{2} \geq 1$ by assumption, and recalling that $\alpha_{3} \geq 1$ and $\beta_{3} \geq 1$. For the coefficient of $(d_1 - m)$ in the last part of inequality \eqref{ineq: F1 case2 genusineq} above, we get
    $$
    (\alpha_{3} + \beta_{3} - 1) d_{1} + (\alpha_{3} - 1) d_{2} \geq d_{1} > 0;
    $$
    and for the coefficient of $(d_2 - 1)$ in the last part of inequality \eqref{ineq: F1 case2 genusineq} above, we get
    $$
    (\alpha_{3} + \beta_{3} - 1) d_{1} + \alpha_{3} d_{2} \geq d_{1} + d_{2} > 0.
    $$
    Therefore, all the equalities above in the inequality \eqref{ineq: F1 case2 genusineq} must hold. We obtain
    $$
    d_{1} = m, 
    \quad 
    d_{2} = 1,
    \quad
    \text{and}
    \quad
     p_a (D^{\nu}) = 0.
    $$

    Substituting these back into the inequality \eqref{ineq: F1 case2 genusineq pre}, we obtain
    \begin{align*}
        0 = p_{a}(D^{\nu}) 
        &\leq \frac{1}{2}(m - 1) m  - \frac{1}{2}(m - 1)((\alpha_{2} + \beta_{2}) m + \alpha_{2} - 1) \\
        &= -\frac{1}{2}(m - 1)((\alpha_{2} + \beta_{2} - 1) m + \alpha_{2} - 1)
        \leq 0.
    \end{align*}

    Therefore, all the equalities here also must hold. If $m > 1$, we must have $\alpha_{2} = 1$ and $\beta_{2} = 0$. If $\{B_1 , B_2\} \subset \lvert C_1 \rvert$, we also need to consider those curves in 
    $\operatorname{Hyp}_{d_1 C_1 + f}(B,2)$ 
    that intersect $B_2$ at two points. To do this, we may switch $B_1$ and $B_2$, then by symmetry, the preceding arguments remain valid. Hence we conclude the proof of part (\ref{prop: F1 Hyp(B,2)=E(B) (3)}) of the proposition.

    Since integral curves in $\lvert C_0 \rvert$, $\lvert C_1 \rvert$ and $\lvert f \rvert$ are already rational, the last statement $\operatorname{Hyp}(B,2) = \mathcal{E}(B)$ 
    immediately follows from part (\ref{prop: F1 Hyp(B,2)=E(B) (2)}) and part (\ref{prop: F1 Hyp(B,2)=E(B) (3)}).
\end{proof}

Now, we will show that $\mathcal{E}(B)$ in Prop.~\ref{prop: F1 Hyp(B,2)=E(B)} is a finite set.

\begin{theorem}\label{thm: F1 bound E(B)}
    Consider $\mathbb{F}_1$. Let $B$ be as in Prop.~\ref{prop: F1 Hyp-C1-(B,2)}. Then:
    \begin{enumerate}[label=(\alph*), ref=(\alph*)]
        \item \label{thm: F1 bound E(B)_a} If $\{B_1 , B_2\} \not\subset \lvert C_{1}\rvert$, then $\mathcal{E}(B)$ is finite and is effectively bounded in terms of $\alpha_{i}$'s and $\beta_{i}$'s.
        \item \label{thm: F1 bound E(B)_b} If $\{B_{1} , B_{2}\} \subset \lvert C_{1}\rvert$, then $\lvert \mathcal{E}(B) \rvert < \infty$.
    \end{enumerate}
\end{theorem}

\begin{proof}
    Combining Prop.~\ref{prop: F1 Hyp-C1-(B,2)} and part (\ref{prop: F1 Hyp(B,2)=E(B) (2)}) of Prop.~\ref{prop: F1 Hyp(B,2)=E(B)}, we see that $\left\lvert \bigcup_{d_1 \geq 1} \operatorname{Hyp}_{d_1 C_{1}}(B,2) \right\rvert$ is finite and effectively bounded in terms of $\alpha_i$'s and $\beta_i$'s.
    
    By Prop.~\ref{prop: F1 Hyp(B,2)=E(B)}, we see it remains to study the cardinality of $\operatorname{Hyp}_{C_{0}}(B,2)$, $\operatorname{Hyp}_{f}(B,2)$, and $\operatorname{Hyp}_{d_{1} C_{1} + f}(B,2)$ with $d_1 \geq 1$.

    Since $\lvert C_0 \rvert =  \{ C_{0} \}$, we have $\lvert \operatorname{Hyp}_{C_{0}}(B,2) \rvert \leq 1$. By the same argument as in the proof of Thm.~\ref{thm: F0 bound E(B)} for the $\mathbb{F}_{0}$ case, we see that any integral curve $D \in \operatorname{Hyp}_{f}(B,2)$ contains some point in $N$. Therefore, we obtain 
    $$
    \lvert \operatorname{Hyp}_{f}(B,2) \rvert \leq \lvert N \rvert 
    = \sum_{i < j} (B_{i} \cdot B_{j}).
    $$

    Now, it remains to study $\operatorname{Hyp}_{d_{1} C_{1} + f}(B,2)$, where $d_{1} \geq 1$. By part (\ref{prop: F1 Hyp(B,2)=E(B) (3)}) of Prop.~\ref{prop: F1 Hyp(B,2)=E(B)}, we know that if 
    $\operatorname{Hyp}_{d_1 C_1 + f}(B,2) \neq \emptyset$ for some $d_1 \geq 1$, then $B_1 \in \lvert C_1 \rvert$. 
    
    For $D \in \operatorname{Hyp}_{d_{1} C_{1} + f}(B,2)$, we first consider the case that $D$ intersects $B_1$ at two points. Then, as in the proof of part (\ref{prop: F1 Hyp(B,2)=E(B) (3)}) of Prop.~\ref{prop: F1 Hyp(B,2)=E(B)}, we have $D \cap B = \{p_{1,2},p_{1,3}\}$ for some points $p_{1,2} \in B_1 \cap B_2$ and $p_{1,3} \in B_1 \cap B_3$, 
    $D \cap B_2 = \{p_{1,2}\}$, 
    $D \cap B_3 = \{p_{1,3}\}$,
    and $D$ intersects $B_1$ transversally at both $p_{1,2}$ and $p_{1,3}$. Moreover, we know that $p_{1,2}$ is a unibranch $d_1$-fold point of $D$, and $B_{2}$ is tangent to $D$ at $p_{1,2}$ since 
    $(B_{2} \cdot D)_{p_{1,2}} 
    = (B_{2} \cdot D) 
    \geq d_1 + 1 
    > d_1$.

    \begin{claim}\label{claim: degree inequality between two unicuspidal curves on F1}
        Fix $p_{1,2}$ and $p_{1,3}$. For any $d_1 \geq 1$, there is at most one integral curve in $\operatorname{Hyp}_{d_1 C_1 + f}(B,2)$ that contains both $p_{1,2}$ and $p_{1,3}$. Moreover, if $D_{1} \in \operatorname{Hyp}_{m_1 C_{1} + f}(B,2)$ and $D_{2} \in \operatorname{Hyp}_{m_2 C_{1} + f}(B,2)$ are two integral curves containing $p_{1,2}$ and $p_{1,3}$ with $m_{1} < m_{2}$, then
        $$
        (\alpha_{3} + \beta_{3}) m_{1} + \alpha_{3} \leq m_{2}.
        $$
    \end{claim}
    \begin{proof}[Proof of the Claim]
        Suppose there are two curves $D, D^{\prime} \in \operatorname{Hyp}_{d_1 C_{1} + f}(B,2)$ containing $p_{1,2}$ and $p_{1,3}$. Then, both of them intersect $B_2$ only at $p_{1,2}$ and intersect $B_3$ only at $p_{1,3}$, $\operatorname{mult}_{p_{1,2}}(D) = \operatorname{mult}_{p_{1,2}}(D^{\prime}) = d_1$, and $p_{1,3}$ is a smooth point on both of them. We obtain
        $$
        (D\cdot B_3)_{p_{1,3}} 
        = (D \cdot B_3) 
        = d_1 (\alpha_3 + \beta_3) + \alpha_3
        = (D^{\prime} \cdot B_3)
        = (D^{\prime} \cdot B_3)_{p_{1,3}},
        $$
        and
        $$
        (D\cdot B_2)_{p_{1,2}} 
        = (D \cdot B_2) 
        = d_1 + 1
        = (D^{\prime} \cdot B_2)
        = (D^{\prime} \cdot B_2)_{p_{1,2}}.
        $$
        By Thm.~\ref{thm: Strong Triangle Inequality}, we obtain
        $$
        \frac{(D\cdot D^{\prime})_{p_{1,2}}}{d_1^2} 
        \geq 
        \frac{(D\cdot B_2)_{p_{1,2}}}{d_1}
        =
        \frac{(D^{\prime} \cdot B_2)_{p_{1,2}}}{d_1}
        =
        \frac{d_1 + 1}{d_1}.
        $$
        Therefore,
        \begin{multline}
            (D \cdot D^{\prime})_{p_{1,3}} 
            \leq (D \cdot D^{\prime}) - (D \cdot D^{\prime})_{p_{1,2}} 
            \leq (d_1^2+2d_1) - (d_1^2 + d_1) \\
            = d_1 
            < d_1 (\alpha_3 + \beta_3) + \alpha_3 
            = (D \cdot B_3)_{p_{1,3}} 
            = (D^{\prime} \cdot B_3)_{p_{1,3}}.
        \end{multline}
        But then, by Thm.~\ref{thm: Strong Triangle Inequality}, the smallest two among
        $$
        (D \cdot D^{\prime})_{p_{1,3}},
        \quad
        (D \cdot B_3)_{p_{1,3}},
        \quad
        (D^{\prime} \cdot B_3)_{p_{1,3}}
        $$
        are equal, contradiction.
        
        Now, let $D_{1} \in \operatorname{Hyp}_{m_1 C_{1} + f}(B,2)$ and $D_{2} \in \operatorname{Hyp}_{m_2 C_{1} + f}(B,2)$ be two curves containing $p_{1,2}$ and $p_{1,3}$ with $m_1 < m_2$. Then both $D_{1}$ and $D_{2}$ are hypertangent to $B_{2}$ at $p_{1,2}$. By Thm.~\ref{thm: Strong Triangle Inequality}, we know that the smallest two among
        $$
        \frac{(D_{1} \cdot D_{2})_{p_{1,2}}}{m_{1}  m_{2}}, \quad \frac{(D_{1} \cdot B_{2})_{p_{1,2}}}{\operatorname{mult}_{p_{1,2}}(D_{1})} = \frac{m_{1} + 1}{m_{1}}, \quad \frac{(D_{2} \cdot B_{2})_{p_{1,2}}}{\operatorname{mult}_{p_{1,2}}(D_{2})} = \frac{m_{2} + 1}{m_{2}}
        $$
        are equal. Hence,
        $$
        \frac{(D_{1} \cdot D_{2})_{p_{1,2}}}{m_{1}  m_{2}} = \frac{m_{2} + 1}{m_{2}}.
        $$

        Moreover, both $D_{1}$ and $D_{2}$ are hypertangent to $B_{3}$ at $p_{1,3}$. Since
        $$
        (D_{1} \cdot D_{2}) = \left( (m_{1} C_{1} + f) \cdot (m_{2} C_{1} + f) \right) = m_{1} m_{2} + m_{1} + m_{2},
        $$
        we have
        $$
        (D_{1} \cdot D_{2})_{p_{1,3}} \leq (D_{1} \cdot D_{2}) - (D_{1} \cdot D_{2})_{p_{1,2}} = (m_{1} m_{2} + m_{1} + m_{2}) - (m_{1} m_{2} + m_{1}) = m_{2}.
        $$

        Furthermore, we have
        $$
        (B_{3} \cdot D_{1})_{p_{1,3}} = (B_{3} \cdot D_{1}) = (\alpha_{3} + \beta_{3}) m_{1} + \alpha_{3}
        \quad
        \text{and}
        \quad
        (B_{3} \cdot D_{2})_{p_{1,3}} = (B_{3} \cdot D_{2}) = (\alpha_{3} + \beta_{3}) m_{2} + \alpha_{3}.
        $$
        
        Then,
        $$
        (B_{3} \cdot D_{1})_{p_{1,3}} = (\alpha_{3} + \beta_{3}) m_{1} + \alpha_{3} < (\alpha_{3} + \beta_{3}) m_{2} + \alpha_{3} = (B_{3} \cdot D_{2})_{p_{1,3}}.
        $$
        Note that $p_{1,3}$ is a smooth point on all three curves $D$, $D^{\prime}$ and $B_3$. By Thm.~\ref{thm: Strong Triangle Inequality}, we must have
        $$
        (\alpha_{3} + \beta_{3}) m_{1} + \alpha_{3} = (B_{3} \cdot D_{1})_{p_{1,3}} = (D_{1} \cdot D_{2})_{p_{1,3}} \leq m_{2}.
        $$
        This completes the proof of the claim.
    \end{proof}

    If $\{B_1 , B_2\} \cap \lvert C_1 \rvert = \{B_1\}$, then, by part (\ref{prop: F1 Hyp(B,2)=E(B) (3)}) of Prop.~\ref{prop: F1 Hyp(B,2)=E(B)}, we obtain 
    $\operatorname{Hyp}_{d_1 C_1 + f} = \emptyset$ for any $d_1 \geq 2$. Moreover, in this case, note that there are at most $(B_1 \cdot B_2)$ choices for $p_{1,2}$ and at most $(B_1 \cdot B_3)$ choices for $p_{1,3}$, and by Claim \ref{claim: degree inequality between two unicuspidal curves on F1}, we know 
    $$
    \lvert \operatorname{Hyp}_{C_1 + f} (B,2) \rvert
    \leq 
    (B_1 \cdot B_2) (B_1 \cdot B_3)
    = (\alpha_2 + \beta_2) (\alpha_3 + \beta_3).
    $$
    Combining Prop.~\ref{prop: F1 Hyp-C1-(B,2)} and Prop.~\ref{prop: F1 Hyp(B,2)=E(B)}, we obtain part \ref{thm: F1 bound E(B)_a} of the theorem.

    If $\{B_1 , B_2\} \subset \lvert C_1 \rvert$, consider the blow-up morphism 
    $\pi\colon \mathbb{F}_{1} \to \mathbb{P}^{2}$ 
    obtained by contracting $C_0$. We have
    $$
    \deg \pi(B_{1}) = \deg \pi(B_{2}) = 1, 
    \quad 
    \deg \pi(B_{3}) = \alpha_{3} + \beta_{3} \geq 2, 
    \quad
    \deg \pi(D) = d_1 + 1,
    $$
    and $\pi(D) \cap \pi(B) = \{q, \pi(p_{1,2}), \pi(p_{1,3})\}$, where $q = \pi(C_0)$ is the blow-up point.

    Then, $\pi(B)$ is still a $3C$-curve. Indeed, in $\mathbb{F}_1$, $B_{1}, B_{2} \in \lvert C_{1} \rvert$ ensures that $B_{1}$ and $B_{2}$ do not intersect $C_{0}$. Now, the Corvaja--Zannier theorem \cite[Thm.~1]{CZ13} tells us
    $$
    d_1 + 1 = \deg \pi(D) \leq \gamma(B),
    $$
    for some constant positive integer $\gamma(B)$ determined by $B$.

    Let $I_1 = \{m_{1}, m_{2}, \dots, m_{n}, \dots\} \subset \mathbb{N}$ be the set of all positive integers such that for each $m_{i} \in I_1$, there exists a curve in $\operatorname{Hyp}_{m_i C_{1} + f}(B,2)$ containing both $p_{1,2}$ and $p_{1,3}$. Then, by Claim \ref{claim: degree inequality between two unicuspidal curves on F1}, for all $i \geq 2$, we obtain
    $$
    m_{i} \geq (\alpha_{3} + \beta_{3}) m_{i-1} + \alpha_{3} > (\alpha_{3} + \beta_{3}) m_{i-1}.
    $$
    By induction, we get
    $$
    m_{i} > (\alpha_{3} + \beta_{3})^{i-1} m_{1}, \quad \forall i \geq 2.
    $$

    Hence,
    $$
    (\alpha_{3} + \beta_{3})^{i-1} m_{1} < \gamma(B) - 1, \quad \forall i \geq 2.
    $$
    If $\lvert I_1 \rvert \geq 2$, this implies
    $$
    \lvert I_1 \rvert = \max_{i \in \mathbb{N}, \ m_i \in I} \{i\} < 1 + \log_{(\alpha_{3} + \beta_{3})} \left( \frac{\gamma(B) - 1}{m_{1}} \right) \leq 1 + \log_{(\alpha_{3} + \beta_{3})} (\gamma(B) - 1).
    $$

    Therefore, we conclude that $I_1$ is a finite set, and 
    $$
    \lvert I_1 \rvert 
    \leq 1 + \lfloor \log_{(\alpha_{3} + \beta_{3})} (\gamma(B) - 1) \rfloor.
    $$

    In the case 
    $\{B_1 , B_2\} \subset \lvert C_1 \rvert$, we also need to consider those curves in  
    $\operatorname{Hyp}_{d_1 C_1 + f}(B,2)$
    which intersect $B_2$ at two points. 
    If $D \in \operatorname{Hyp}_{d_1 C_1 + f}(B,2)$ is such a curve, then $D \cap B_2 = \{p_{1,2} , p_{2,3}\}$,
    where $\{p_{1,2}\} = B_1 \cap B_2$ 
    and $p_{2,3} \in B_2 \cap B_3$. 
    By symmetry, after switching $B_1$ and $B_2$, the preceding arguments remain valid. Let $I_2 = \{m_{1}^{\prime}, m_{2}^{\prime}, \dots, m_{n}^{\prime}, \dots\} \subset \mathbb{N}$ be the set of all positive integers such that for each $m_{i}^{\prime} \in I_2$, there exists a curve in $\operatorname{Hyp}_{m_i^{\prime} C_{1} + f}(B,2)$ containing both $p_{1,2}$ and $p_{2,3}$, then, by symmetry, we also have
    $$
    \lvert I_2 \rvert 
    \leq 1 + \lfloor \log_{(\alpha_{3} + \beta_{3})} (\gamma(B) - 1) \rfloor.
    $$

    Now, since $\{p_{1,2}\} = B_1 \cap B_2$ and there are only $\alpha_3 + \beta_3$ choices for $p_{1,3}$ and $\alpha_3 + \beta_3$ choices for $p_{2,3}$, combining Prop.~\ref{prop: F1 Hyp(B,2)=E(B)},  we conclude that $\mathcal{E}(B)$ is a finite set, which gives part \ref{thm: F1 bound E(B)_b} of the theorem and completes the proof.
\end{proof}

\begin{example}
    In the case where $B_{1}, B_{2} \in \lvert C_{1} \rvert$, if $B_{3} \in \lvert C_{1} + f \rvert$ and $B$ is in general position, with all the notations being the same as above in the proof of Thm.~\ref{thm: F1 bound E(B)}, then, by another result of Corvaja--Zannier \cite[Thm.~1.1]{CZ08}, we have $\gamma(B) \leq 2^{15} \cdot 35$. Then, we obtain
    $$
    \lvert I_{j} \rvert 
    < 1 + \log_{2} \left( 2^{15} \cdot 35 - 1 \right) 
    < 1 + \log_{2} \left( 2^{15} \cdot 35 \right) 
    = 1 + 15 + \log_{2}(35) < 22,
    \quad
    j = 1,2.
    $$
\end{example}

\begin{proposition}\label{prop: emptyness of exceptional set F1}
    Consider $\mathbb{F}_1$. Let $B$ be as in Prop.~\ref{prop: F1 Hyp-C1-(B,2)}. 
    
    If $\alpha_i \geq 3$ for all $i \in \{1,2,3\}$, $\beta_l \geq 3$ for some $l \in \{1,2,3\}$, and $B$ is general, then
    $
    \mathcal{E}(B) = \emptyset
    $.
\end{proposition}
\begin{proof}
    The statement follows by combining Prop.~\ref{prop: F1 Hyp(B,2)=E(B)} and Thm.~\ref{thm: F1 bound E(B)} with the argument from the proof of Prop.~\ref{prop: emptyness of exceptional set Fe}.
\end{proof}

Finally, we provide an example that if the assumption on the $3C$-curve $B$---that none of its three components lies in $\lvert f \rvert$ or is $C_{0}$---is not satisfied, then $\operatorname{Hyp}(B,2)\neq \mathcal{E}(B)$:
\begin{example}\label{exm:failure-of-hyper-bitangency-implies-rationality-on-F1-prime}
    Let $D \subset \mathbb{P}^2$ be a smooth cubic curve with a flex fixed as the neutral element for the group structure on it, and let $p \in D$ be a $9$-torsion point which is not a flex. By \cite[Thm.~1.3, Prop.~2.2]{Moschetti2025}, there exists a unique pencil $ \mathcal{V} $ containing $D$ such that $p$ is the only base point of the pencil.
        
    Let $q$ be a flex of $D$, and let $L$ be the tangent line to $D$ at $q$.
    \begin{claim}
        There exists a Zariski-open dense subset $U \subset \mathcal{V}$ where every member intersects $L$ transversally. 
    \end{claim}
    
    Let $C$ be a smooth cubic curve in $U$; then $C\cap D = \{p\}$ by construction.  Blow up a point $x$ in $L \setminus (C \cup D)$. In the blow-up surface $\operatorname{Bl}_x\mathbb{P}^2 \cong \mathbb{F}_{1}$, denote the exceptional line by $E$. The curve $B \coloneqq \tilde{L} \cup E \cup \tilde{C}$ is a $3C$-curve, and by Prop.~\ref{prop: Criterion for bigness on Hirzebruch Surfaces} we know 
    $$
    K_{\mathbb{F}_1}+\mathcal{L}(B) 
    \sim -2C_1 - f+f + C_0 + 3C_1 
    \sim C_1+C_0
    \sim 2C_1 - f
    $$ 
    is a big divisor on $\mathbb{F}_{1}$. However, $\tilde{D} \in \operatorname{Hyp}(B,2)$ has geometric genus $1$. Therefore, $\operatorname{Hyp}(B,2) \neq \mathcal{E}(B)$ in this case.
    \begin{figure}[htbp]
        \centering
        \includegraphics[width=0.45\textwidth]{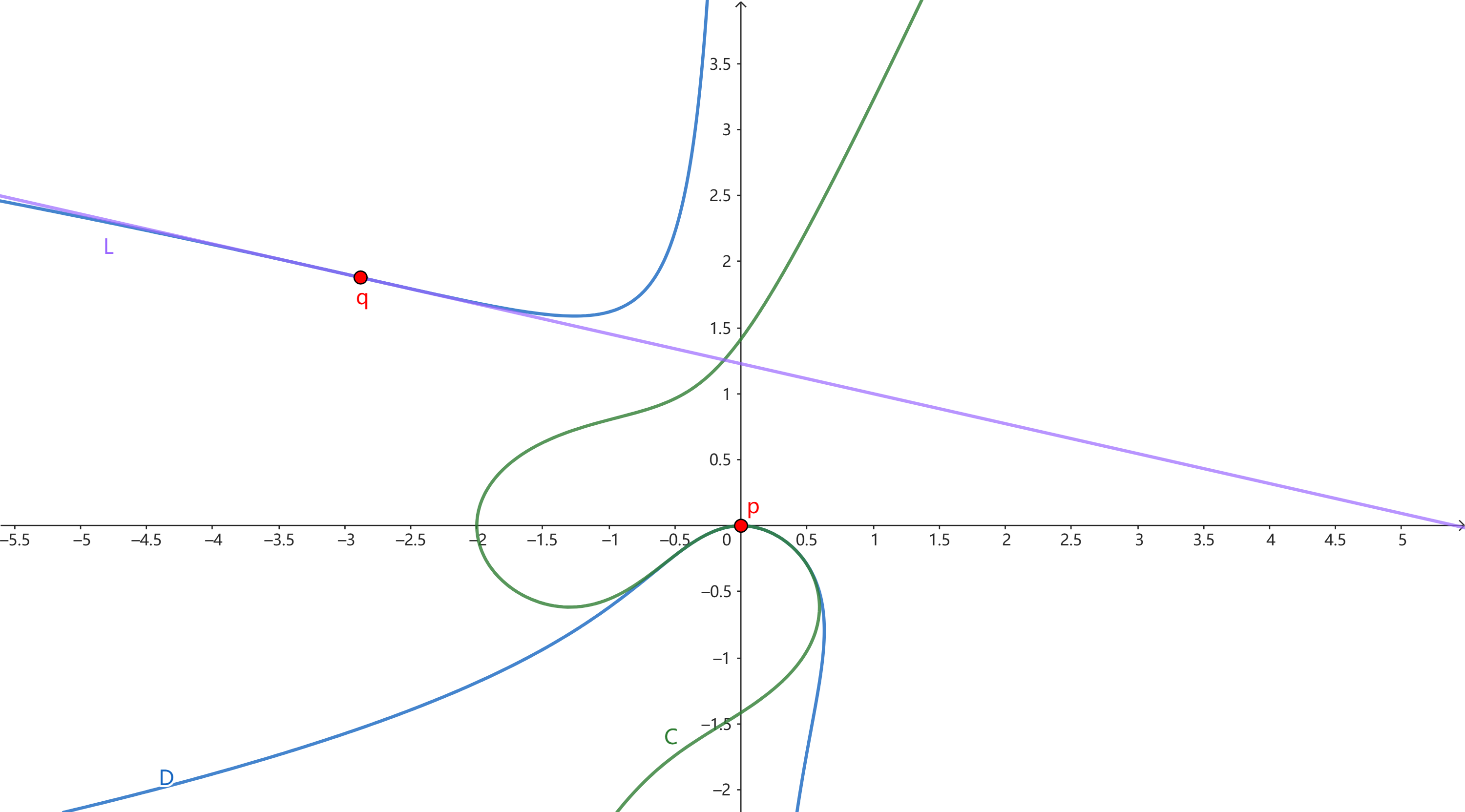}
        \caption{Example \ref{exm:failure-of-hyper-bitangency-implies-rationality-on-F1-prime}}
        \label{fig:counterexample}
    \end{figure}
    
    It remains to prove the claim. By construction, the pencil $\mathcal{V}$ has only one base point $p$ which is not contained in the line $L$. Therefore, the restriction of $\mathcal{V}$ to $L$---denoted as $\mathcal{V}_{L}$---is a base-point-free linear system on $L$. By Bertini's theorem (\cite[Ch.~\Romannum{3}, Cor.~10.9, p.~274]{Hartshorne1977}), there exists a Zariski-open dense subset $U^{\prime} \subseteq \mathcal{V}_{L}$ where every member is smooth as a closed subscheme of $L$. This implies that there exists a Zariski-open dense subset $U \subseteq \mathcal{V}$ where every member intersects $L$ transversally. Thus, the claim is established.
\end{example}

\bibliographystyle{plain} 
\bibliography{ref}

\end{document}